\documentclass[11pt,reqno,a4paper]{amsart}
\usepackage[T1]{fontenc}
\usepackage{lmodern}
\usepackage{amsmath,amssymb,mathtools,mathrsfs}
\usepackage{microtype,booktabs,longtable,needspace}
\usepackage{tikz}
\usetikzlibrary{arrows.meta}
\usepackage[margin=1.1in]{geometry}
\usepackage{xcolor}
\usepackage{hyperref}
\hypersetup{colorlinks,	linkcolor={red!50!black},	citecolor={blue!50!black},	urlcolor={blue!80!black},
 pdftitle={Hessian Defect, Compatibility Degree, and Canonical Decomposition},pdfauthor={Jiarui Fei}}

\calclayout
\allowdisplaybreaks[2]

\newtheorem{theorem}{Theorem}[section]
\newtheorem{lemma}[theorem]{Lemma}
\newtheorem{proposition}[theorem]{Proposition}
\newtheorem{corollary}[theorem]{Corollary}
\theoremstyle{definition}
\newtheorem{definition}[theorem]{Definition}
\newtheorem{example}[theorem]{Example}
\theoremstyle{remark}
\newtheorem{remark}[theorem]{Remark}
\newtheorem{problem}[theorem]{Problem}
\newtheorem*{maintheoremA}{Theorem A}
\newtheorem*{maintheoremB}{Theorem B}
\newtheorem*{maintheoremC}{Theorem C}
\numberwithin{equation}{section}

\DeclareMathOperator{\Hom}{Hom}
\DeclareMathOperator{\End}{End}
\DeclareMathOperator{\im}{im}
\DeclareMathOperator{\coker}{coker}
\DeclareMathOperator{\rad}{rad}
\DeclareMathOperator{\rank}{rank}
\DeclareMathOperator{\grk}{grk}
\DeclareMathOperator{\PHom}{PHom}
\DeclareMathOperator{\IHom}{IHom}
\DeclareMathOperator{\PC}{PC}
\DeclareMathOperator{\Gr}{Gr}
\DeclareMathOperator{\rep}{rep}
\DeclareMathOperator{\gen}{gen}
\DeclareMathOperator{\ext}{ext}
\DeclareMathOperator{\tors}{tors}

\newcommand{\op}[1]{\operatorname{#1}}
\newcommand{\kk}{\mathbb C}
\newcommand{\Z}{\mathbb Z}
\newcommand{\E}{\op{E}}
\newcommand{\ec}{\check{\mathsf e}}
\newcommand{\ee}{\mathsf e}
\newcommand{\dv}{\underline{\dim}}
\newcommand{\den}{\underline{\operatorname{den}}}
\newcommand{\dT}{\underline{\operatorname{den}}^{\mathrm{trop}}}
\newcommand{\rh}{\underline{\operatorname{rank}}}
\newcommand{\fc}{\check f}
\newcommand{\dtc}{{\check\delta}}
\newcommand{\etc}{{\check\eta}}

\newcommand{\betac}{\check\beta}
\newcommand{\M}{\mathcal M}
\newcommand{\N}{\mathcal N}
\newcommand{\HH}{\mathscr H}
\newcommand{\symE}{\mathsf e_{\mathrm{sym}}}
\newcommand{\inn}{\mathrm{in}}
\newcommand{\out}{\mathrm{out}}
\newcommand{\nf}{\mathrm{nf}}
\newcommand{\tr}{\mathsf T}
\newcommand{\CC}{C_{\gen}}
\newcommand{\Uca}{\mathcal U}
\newcommand{\bx}{\boldsymbol x}
\newcommand{\by}{\boldsymbol y}
\newcommand{\pair}[2]{\langle #1,#2\rangle}
\newcommand{\comp}[2]{(#1\parallel #2)}
\newcommand{\mneg}{m^-}
\newcommand{\arxiv}[1]{\href{https://arxiv.org/abs/#1}{arXiv:#1}}

\begin{document}
\title[Hessian defect and canonical decomposition]{Hessian Defect, Compatibility Degree, and Canonical Decomposition}
\author{Jiarui Fei}
\address{School of Mathematical Sciences, Shanghai Jiao Tong University}
\email{jiarui@sjtu.edu.cn}
\thanks{The author was supported in part by the National Natural Science Foundation of China (Nos.~12131015 and 12571038).}
\date{}

\begin{abstract}
We express the difference between dimension and denominator vectors as the rank defect of a Hessian differential. For a very general Jacobi-finite potential, this defect vanishes on a general representation in a principal component precisely when the positive-length cycles at the vertex act trivially. The associated Hessian vectors satisfy tropical $X$-mutation and Auslander--Reiten translation. We also give two distinct cluster variables with the same denominator vector but different Hessian vectors. Their Hom pairing extends ordered denominator compatibility, and its negative part computes the canonical multiplicity of any extended-reachable indecomposable class in an arbitrary presentation weight.
\end{abstract}

\maketitle
\setcounter{tocdepth}{2}
\tableofcontents
\enlargethispage{2pt}

\section*{Introduction}

\subsection{Denominator vectors and Hessian defects}
The denominator vector need not equal the dimensions of the corresponding representation. Fu--Keller prove the denominator--dimension inequality in a $2$-Calabi--Yau setting \cite[Proposition 6.8]{FK10}; Derksen--Weyman--Zelevinsky prove the QP form used here and give a necessary subrepresentation condition for equality \cite[Corollary 5.5]{DWZ2}. Equality has also been studied for acyclic cluster algebras \cite{BMR09} and for cluster algebras from triangulated surfaces \cite{Yur24}. We study this difference by comparing the rank of the Hessian differential with the general rank in the same presentation space.

Let $J=J(Q,W)$ be a finite-dimensional Jacobian algebra of a nondegenerate quiver with potential, and let $b_i$ be the $i$-th row of its exchange matrix. The projective resolution of the simple at $i$ contains the Hessian presentation
\[ p_i:P_{\inn}^i\longrightarrow P_{\out}^i,\qquad \coker p_i=\rad P_i. \]
Its projective terms are recalled in Section~\ref{sec:QP}. For a module $N$, write $p_N^*$ for precomposition by $p$, and define the \emph{local Hessian defect} by
\[ h_i(N)=\max_{p\in\Hom_J(P_{\inn}^i,P_{\out}^i)}\rank p_N^*-\rank(p_i)_N^*. \]
Here $\hom(b_i,N)$ denotes the generic value of $\dim_\kk\Hom_J(\coker p,N)$ in this presentation space. General presentations and their principal components are recalled in Section~\ref{sec:general}.

\begin{maintheoremA}[Theorems~\ref{thm:local}, \ref{thm:local-radical-bound}, \ref{thm:cycle-action}, and~\ref{thm:denominator}]
For every finite-dimensional $J$-module $N$,
\[ h_i(N)=\hom(\rad P_i,N)-\hom(b_i,N),\qquad \max_{c\in\rad(e_iJe_i)}\rank R_c\leq h_i(N), \]
where $R_c$ is right multiplication by $c$ on $N(i)$. There is a nonempty very general set $\mathcal V_Q$ of Jacobi-finite nondegenerate potentials such that, for $W\in\mathcal V_Q$ and the general ordinary part $N$ of every principal component,
\[ 0\leq h_i(N)\leq\dim_\kk\bigl(N(i)\rad(e_iJe_i)\bigr),\qquad h_i(N)=0\Longleftrightarrow N(i)\rad(e_iJe_i)=0. \]
These assertions hold in every mutation-equivalent QP. If $M$ corresponds to a noninitial cluster variable $z$ in the seed $t$, then
\[ \den_t(z)(i)=\dim M(i)-h_i(M). \]
\end{maintheoremA}

For ordinary reachable representations, Corollary~\ref{cor:regularity} also identifies zero defect with the subrepresentation condition of \cite[Corollary 5.5]{DWZ2}. For the remaining main statements, we fix one potential $W\in\mathcal V_Q$; the same choice applies in every mutation-equivalent QP.

The Hessian-vector construction, its extended mutation, and the canonical-multiplicity formula use generic Hom dimensions and do not assume the generic-pairing conjecture. To identify the local Hom expression with a tropical denominator, one needs $f_N(b_i)=\hom(b_i,N)$, where $f_N$ is the tropical $F$-polynomial of Section~\ref{ss:tropical-F}. The corresponding Laurent identity also requires that the terms of minimal exponent do not cancel after monomial specialization. Both conditions hold for reachable cluster monomials; for general characters they remain explicit in Theorem~\ref{thm:denominator}.

\subsection{Hessian vectors}
The exponent of $x$ in the denominator of $z$ defines an ordered compatibility degree $d(x,z)$, independent of the other variables in a cluster containing $x$ \cite{CL18,CL}. We use the signed convention $d(x,x)=-1$. In finite type this recovers the compatibility degree of Fomin--Zelevinsky \cite{FZ032,CP15}; related descriptions of cluster complexes use Coxeter-sortable elements \cite{Rea07} and combinatorial frameworks \cite{RSpeyer}. The Hessian vector will recover the same degree with its arguments reversed.

For a presentation weight $\eta$, let $\dv(\eta)$ be its general ordinary dimension vector and $\mneg(\eta)$ the vector of canonical multiplicities of the negative simples, as in Section~\ref{sec:general}. An indecomposable class is \emph{extended-reachable} if a sequence of ordinary mutations and decorated Auslander--Reiten translations takes it to a negative simple. Choose a \emph{negative normalization} $(w,i)$, so that $w$ takes $\eta$ to $-\mathbf e_i$, and write $w(S_j^-)=(N_j,N_j^-)$ over the resulting algebra $J^w$. For each $j$, we evaluate the defect of $N_j$ at $i$ and set
\[ \Delta^w(\eta)(j)=h_i^{J^w}(N_j),\qquad \gamma^w(\eta)=\dv(\eta)-\mneg(\eta)-\Delta^w(\eta). \]

\begin{maintheoremB}[Theorems~\ref{thm:endpoint-rank}, \ref{thm:intrinsic-vector}, and~\ref{thm:ordinary-vector}]
The vectors $\Delta^w(\eta)$ and $\gamma^w(\eta)$ are independent of the negative normalization. Write them as $\Delta(\eta)$ and $\gamma(\eta)$. At a negative simple, $\gamma(-\mathbf e_i)=-\mathbf e_i$. Ordinary mutation at $k$ leaves all other coordinates unchanged and gives
\[ \gamma'(k)=\max\left\{\sum_j[b_{kj}]_+\gamma(j),\sum_j[-b_{kj}]_+\gamma(j)\right\}-\gamma(k). \]
For the action $\tau_y$ of Auslander--Reiten translation on weights,
\[ \gamma(\tau_y\eta)=\dv(\gamma(\eta)B)-\gamma(\eta), \]
where $\gamma(\eta)B$ is a coweight. The assignment $\eta\mapsto\gamma(\eta)$ is injective. If $\eta$ has ordinary part $M\ne0$, then
\[ 0\leq\gamma(\eta)\leq\dv M,\qquad \gamma(\eta)(j)>0\Longleftrightarrow M(j)\ne0. \]
For a reachable $\eta$, with cluster variable $z=X(\eta)$,
\[ \gamma(\eta)(j)=d(z,x_{j;t}),\qquad \den_t(z)(j)=d(x_{j;t},z). \]
\end{maintheoremB}

Thus the two vectors use opposite orders of the same compatibility degree. In matrix form, the relation is transposition together with interchange of the two seeds, not symmetry of $d$; see Section~\ref{ss:ordinary-denominator}.

The distinction is visible in rank four. In Example~\ref{ex:equal-denominators}, we give distinct cluster variables $z,z'$ with
\[ \den_t(z)=\den_t(z')=(4,6,4,6), \]
but Hessian vectors $(4,6,4,5)$ and $(4,5,4,6)$. This gives a counterexample to the denominator-vector conjecture of Fomin--Zelevinsky \cite[Conjecture 4.17]{FZ031}.

\subsection{Compatibility and canonical decomposition}
Let $m_\eta(\delta)$ be the multiplicity of $\eta$ in the Derksen--Fei canonical decomposition of a presentation weight $\delta$ \cite{DF}. For an extended-reachable indecomposable $\eta$ and arbitrary $\delta$, define
\[ \comp{\eta}{\delta}=\mathcal I_B(\gamma(\eta),\delta)=\hom(\delta,\gamma(\eta)B)-\pair{\gamma(\eta)}{\delta}. \]
The second Hom argument is a coweight. The generic Hom and extension numbers are those of Sections~\ref{sec:QP}--\ref{sec:general}; set $\symE(\eta,\delta)=\ee(\eta,\delta)+\ee(\delta,\eta)$.

\begin{maintheoremC}[Theorems~\ref{thm:pair-identity} and~\ref{thm:multiplicity}]
The degree $\comp{\eta}{\delta}$ is invariant under simultaneous extended mutation and additive on canonical decompositions in $\delta$. There is a nonnegative defect $h(\eta,\delta)$ such that
\[ \comp{\eta}{\delta}=\symE(\eta,\delta)-m_\eta(\delta)-h(\eta,\delta), \]
and
\[ m_\eta(\delta)=[-\comp{\eta}{\delta}]_+. \]
Moreover,
\[ \comp{\eta}{\delta}\leq0\Longleftrightarrow\symE(\eta,\delta)=0,\qquad \comp{\eta}{\delta}\geq0\Longleftrightarrow m_\eta(\delta)=0. \]
For ordinary cluster variables, $\comp{\eta}{\delta}=d(X(\eta),X(\delta))$.
\end{maintheoremC}

To prove the multiplicity formula, we apply a negative normalization to both classes. At the final vertex, the degree is positive if the transformed second class has nonzero ordinary support; otherwise it is the negative of its decoration multiplicity. The sequence preserves canonical multiplicities, so the local sign theorem gives the formula for arbitrary $\delta$.

Section~\ref{sec:compat} also compares this degree with $f$-compatibility \cite{FuGyoda} and the Hessian vector with Schur reductions \cite{FeiProjection}. The latter comparison gives $\gamma(M)=\dv M$ for extended-reachable Schur representations and yields the acyclic multiplicity formula in Corollary~\ref{cor:acyclic-multiplicity}.

\subsection{Conventions and notation}
All modules are right modules over $\kk$, all vectors are rows, and
\[ b_{ij}=\#\{i\to j\}-\#\{j\to i\},\qquad \pair{u}{v}=uv^{\tr}. \]
The symbol $e_i$ is a vertex idempotent and $\mathbf e_i$ is a standard basis row. Write $[a]_+=\max\{a,0\}$, coordinatewise for vectors. The symbols $\hom$ and $\ee$ denote Hom and presentation $\E$-space dimensions, and $\dv M$ is the dimension vector of a module. We use matching letters $\delta,\dtc$ for the weight and coweight of the same representation or principal component. In a pairing with independently variable arguments, we write $\delta$ for the source weight and $\etc$ for the target coweight. The notation $\dv(\delta)$, respectively $\dv(\dtc)$, denotes the generic ordinary dimension vector of the principal component labelled by the weight $\delta$, respectively the coweight $\dtc$. When $\delta$ and $\check\delta$ label the same principal component, they are related by \eqref{eq:weight-coweight}. The row $b_i$ is a weight and $-b_i$ a coweight. Ordinary and negative parts of a decorated representation are kept separate. For a module map $H$, $\rh H=(\rank H(j))_j$ is its dimension-vector rank.

The subscripts on an invariant specify its algebra or cluster seed when necessary. Mutation sequences written as lists are applied from left to right. A general value is taken on a nonempty open subset of the indicated irreducible parameter space. In the main text, the standing choice $W\in\mathcal V_Q$ begins after Theorem~\ref{thm:simultaneous-zero-support}.

\begingroup\small
\setlength{\LTpre}{5pt}\setlength{\LTpost}{5pt}
\renewcommand{\arraystretch}{1.12}
\begin{longtable}{@{}p{4.1cm}p{\dimexpr\linewidth-4.1cm-2\tabcolsep\relax}@{}}
\toprule Symbol & Meaning\\\midrule
\endfirsthead
\toprule Symbol & Meaning\\\midrule
\endhead
\bottomrule\endfoot
$\PHom(\delta),\IHom(\dtc)$ & Projective and injective presentation spaces.\\
$\PC(\delta)$ & Principal component of the weight $\delta$.\\
$\dv(\delta),\dv(\dtc)$ & Generic ordinary dimension vectors of a weight and a coweight.\\
$\mneg_i(\delta),m_\eta(\delta)$ & Canonical multiplicities of $-\mathbf e_i$ and $\eta$.\\
$d_N^*,\Hom(d,N),\E(d,N)$ & Precomposition map, its kernel, and its cokernel.\\
$\ee,\ec,\symE$ & Dimension of $\E$-invariants and $\symE(\eta,\delta)=\ee(\eta,\delta)+\ee(\delta,\eta)$.\\
$p_i,\Gamma_i(N)$ & Hessian presentation and its evaluation on $N$.\\
$\HH_i(N),r_i(N),h_i(N)$ & Evaluated presentation space, its general rank, and local Hessian defect.\\
$\grk_J(X,Y),\rh H$ & General dimension-vector rank in $\Hom_J(X,Y)$, and rank vector of $H$.\\
$q_i(\delta)$ & $\hom(\delta,-b_i)+\delta(i)$; the second argument is a coweight.\\
$A_i,\mathfrak r_i,\kappa_i(N)$ & $A_i=e_iJe_i$, its radical, and $\dim_\kk(N(i)/N(i)\mathfrak r_i)$.\\
$\boldsymbol h(\delta)$ & Vector of general local Hessian defects.\\
$\boldsymbol\kappa(M)$ & The row vector $(\kappa_i(M))_{i\in Q_0}$.\\
$\rho_+(M),\rho_-(M)$ & Quotient by the images, and intersection of the kernels, of radical endomorphisms.\\
$\gamma(\eta),\Delta(\eta)$ & Hessian vector and defect vector of an extended-reachable indecomposable.\\
$\mu_k^y,\mu_k^{\check y},\mu_k^x$ & Mutation of weights, coweights, and tropical $X$-vectors.\\
$\tau_y,\tau_{\check y},\tau_x$ & Corresponding Auslander--Reiten operations.\\
$(w,i),w_y,w_x$ & Negative normalization and the two actions of its sequence.\\
$\mathcal I_B(\gamma,\delta)$ & $\hom(\delta,\gamma B)-\pair{\gamma}{\delta}$.\\
$(\eta\parallel\delta),h(\eta,\delta)$ & Signed compatibility degree and its scalar Hessian defect.\\
$c_f(\eta,\delta)$ & Fu--Gyoda $f$-compatibility degree of ordinary cluster variables.\\
$d(x,z),\den_t,\dT$ & Ordered denominator degree, Laurent denominator vector, and denominator before monomial specialization.\\
$\mathcal V_Q$ & The very general initial-potential set.\\
\end{longtable}
\endgroup

Sections~\ref{sec:QP}--\ref{sec:cluster} recall the background. Section~\ref{sec:den} treats local defects, Section~\ref{sec:vectors} constructs Hessian vectors, and Section~\ref{sec:projective} gives projective and injective realizations. The compatibility and multiplicity results are in Section~\ref{sec:compat}, followed by further questions in Section~\ref{sec:outlook}.

\section{Preliminaries on quivers with potentials}\label{sec:QP}
\subsection{The local Jacobian complex}
Let $(Q,W)$ be a reduced nondegenerate QP and let $J=J(Q,W)$ be its finite-dimensional completed Jacobian algebra. Write $P_i=e_iJ$, $I_i=D(Je_i)$, and $S_i=P_i/\rad P_i$, where $D=\Hom_\kk(-,\kk)$. For a nonnegative vector $\beta$, write $P(\beta)=\bigoplus_iP_i^{\oplus\beta(i)}$ and similarly for $I(\beta)$. The reduced quiver has no loops or oriented $2$-cycles. We use the mutation theory of Derksen--Weyman--Zelevinsky \cite{DWZ1,DWZ2}.

We put
\[ P_{\inn}^i=\bigoplus_{a:j\to i}P_j,\qquad P_{\out}^i=\bigoplus_{b:i\to k}P_k. \]
The cyclic second derivatives of \cite[(4.4)--(4.7)]{DWZ2} give the Hessian presentation
\[ p_i=(\partial_{ab}W)_{b,a}:P_{\inn}^i\longrightarrow P_{\out}^i, \]
where an entry corresponding to $a:j\to i$ and $b:i\to k$ is a path combination from $k$ to $j$, viewed as a map $P_j\to P_k$.
There is an exact sequence
\begin{equation} \label{eq:simple-resolution} P_{\inn}^i\xrightarrow{p_i}P_{\out}^i \longrightarrow P_i\longrightarrow S_i\longrightarrow0. \end{equation}
In particular, $\coker p_i=\rad P_i$.

\subsection{Decorated representations and Betti vectors}
A \emph{decorated representation} is a pair $\N=(N,N^-)$, where $N\in\rep J$ and $N^-$ is a finite-dimensional $\kk^{Q_0}$-module. The negative simple $S_i^-$ is $(0,\kk e_i)$. Set
\[ N_{\inn}^i=\bigoplus_{a:j\to i}N(j),\qquad N_{\out}^i=\bigoplus_{b:i\to k}N(k). \]
Use the DWZ maps $\alpha_i:N_{\inn}^i\to N(i)$, $\beta_i:N(i)\to N_{\out}^i$, and $\Gamma_i=(p_i)_N^*$ from \cite[(1.8)--(1.12)]{DWZ2}. The Betti vectors and their duals are
\begin{align*}
\beta_{+,\N}(i)&=\dim\coker\alpha_i,&
 \beta_{-,\N}(i)&=\dim(\ker\alpha_i/\im\Gamma_i)+\dim N^-(i),\\
 \betac_{+,\N}(i)&=\dim\ker\beta_i,&
 \betac_{-,\N}(i)&=\dim(\ker\Gamma_i/\im\beta_i)+\dim N^-(i).
\end{align*}
The weight and coweight are $\delta_\N=\beta_{+,\N}-\beta_{-,\N}$ and $\dtc_\N=\betac_{+,\N}-\betac_{-,\N}$. Subtraction gives
\begin{equation} \label{eq:coindex-local} \dtc_\N(i)=\dim N(i)-\dim N_{\out}^i+\rank\Gamma_i-\dim N^-(i). \end{equation}
Likewise $\delta_\N(i)=\dim N(i)-\dim N_{\inn}^i+\rank\Gamma_i-\dim N^-(i)$. Hence
\begin{equation} \label{eq:weight-coweight} \dtc_\N=\delta_\N+\dv N\,B. \end{equation}
If primes denote mutation at $k$, the ordinary Hom dimensions satisfy \cite[Lemma 2.6(1)]{FeiRank}
\[ \hom(M',N')-\hom(M,N) =\beta_{-,\M}(k)\betac_{-,\N}(k) -\beta_{+,\M}(k)\betac_{+,\N}(k). \]

\subsection{Presentation spaces}
For a projective presentation $d:P_-\to P_+$, define $\Hom(d,N)$ and $\E(d,N)$ by
\begin{equation} \label{eq:presentation-hom} 0\to\Hom(d,N)\to\Hom(P_+,N)\xrightarrow{d_N^*} \Hom(P_-,N)\to\E(d,N)\to0. \end{equation}
The map $d_N^*$ is precomposition with $d$. For an injective presentation $\check d:I_+\to I_-$, define $\Hom(M,\check d)$ and $\check\E(M,\check d)$ as the kernel and cokernel of $\Hom(M,I_+)\to\Hom(M,I_-)$. For $\M=(M,M^-)$, let $d_\M$ and $\check d_\M$ be the minimal projective and injective presentations, including the summands $P(\dv M^-)\to0$ and $0\to I(\dv M^-)$, respectively. Set $\E(\M,\N)=\E(d_\M,N)$ and $\check\E(\M,\N)=\check\E(M,\check d_\N)$; for presentations, $\E(d,d')=\E(d,\coker d')$. These are the presentation extension spaces of \cite[Section 3]{DF}. Their dimensions satisfy
\[ \ee(\M,\N)=\hom(M,N)-\pair{\delta_\M}{\dv N},\qquad \ec(\M,\N)=\hom(M,N)-\pair{\dtc_\N}{\dv M}. \]
The symmetric invariant is
$\symE(\M,\N)=\ee(\M,\N)+\ee(\N,\M)$.

\subsection{Auslander--Reiten translation}
Write $\nu=D\Hom_J(-,J)$. For $\M=(M,V)$, decompose
$M=M_0\oplus P(\pi)$, where $M_0$ has no projective summand, and set
\[ \tau\M=(\tau M_0\oplus I(\dv V),\kk^{\pi}). \]
Here $\kk^{\pi}$ denotes the graded decoration of dimension $\pi$, and $\tau M_0$ is the kernel of the Nakayama transform of a minimal projective presentation. The inverse uses minimal injective presentations and $\nu^{-1}$. In particular,
\[ \tau P_i=S_i^-,\qquad \tau S_i^-=I_i. \]
The operations are inverse bijections on isomorphism classes, respect direct sums, and commute with ordinary mutation \cite[Proposition 7.10]{DF}.

The presentation version of Auslander--Reiten duality identifies (see also \cite{DWZ2,DF})
\[ \E(\M,\N)\simeq D\Hom(N,\widehat\tau\M),\qquad \check\E(\M,\N)\simeq D\Hom(\widehat\tau^{-1}\N,M). \]
Here, the hat on $\tau$ means that we forget the decorated part of $\tau\M$, and $\Hom(\M,\N)=\Hom_J(M,N)$.

We will use the following invariance of the symmetric extension number.

\begin{lemma}[{\cite[Theorem 7.1]{DWZ2}, \cite[Corollary 7.6]{DF}}]\label{lem:symmetric-invariance}
For decorated representations $\M,\N$, the number $\symE(\M,\N)$ is invariant under simultaneous ordinary mutation and under simultaneous $\tau^{\pm1}$.
\end{lemma}

In particular, representations obtained from negative simples by these operations are rigid.

\section{General presentations and tropical \texorpdfstring{$F$}{F}-polynomials}\label{sec:general}
\subsection{General values and coherence}\label{ss:principal-components}
For a weight $\delta$, set
\[ P_-(\delta)=P([-\delta]_+),\qquad P_+(\delta)=P([\delta]_+), \qquad\PHom(\delta)=\Hom(P_-(\delta),P_+(\delta)). \]
For a coweight $\dtc$, set
\[ \IHom(\dtc)=\Hom(I([\dtc]_+),I([-\dtc]_+)). \]
On a dense open subset of the presentation space, the ordinary cokernel dimension and the negative decoration are constant. The closure of the corresponding isomorphism classes is the \emph{principal component}, denoted by $\PC(\delta)$. It is an irreducible component of the corresponding decorated representation variety: forgetting the fixed decoration identifies it with an irreducible component of $\rep_{\dv(\delta)}(J)$. This follows from the openness of the cokernel incidence projection \cite[Theorem 2.3]{DF}; see also \cite[Definition 3.5 and Theorem 3.6]{FeiRank}. We call its general ordinary part a general representation of weight $\delta$. A coweight labels the same family by general injective presentations.

For a weight $\delta$ and a coweight $\dtc$, write
\[ \dv(\delta)=\dv\coker(\delta),\qquad \dv(\dtc)=\dv\ker(\dtc) \]
for the generic ordinary dimension vectors of their principal components. General projective and injective presentations labelling the same principal component give the same dimension vector.

For fixed $N$, a general $d\in\PHom(\delta)$ gives constant values $\hom(\delta,N)$ and $\ee(\delta,N)$. The same convention is used for a general injective presentation. If both arguments vary, ``general'' means a general pair in the product of the two parameter spaces. Even when the two labels correspond to the same principal component, this is not a diagonal pair $(M,M)$; the generic Hom dimension need not be the general endomorphism dimension.

\begin{lemma}[{\cite{DF}}]\label{lem:cancellation}
For projective multiplicity vectors $\beta_-,\beta_+$, a general map in
$\Hom(P(\beta_-),P(\beta_+))$ is homotopy equivalent to a general presentation in
$\PHom(\beta_+-\beta_-)$. In particular their generic ordinary cokernel dimensions and presentation Hom dimensions agree.
\end{lemma}

After removing contractible summands, a general presentation has disjoint positive and negative projective multiplicities. General decorated representations are both weight- and coweight-coherent:
\[ \beta_\pm=[\pm\delta]_+,\qquad \betac_\pm=[\pm\dtc]_+. \]
The weight and coweight labelling the same principal component are related by \eqref{eq:weight-coweight}. Generality is preserved by mutation \cite[Theorem 1.10]{GLFS} and by the weight--coweight correspondence \cite[Theorem 3.11]{FeiRank}. Consequently, for an independently chosen source weight $\delta$ and target coweight $\etc$, the generic Hom-mutation formula is
\begin{equation} \label{eq:generic-Hom-mutation} \hom(\delta',\etc')-\hom(\delta,\etc) =[-\delta(k)]_+[-\etc(k)]_+ -[\delta(k)]_+[\etc(k)]_+. \end{equation}

The generic weight and coweight mutations are
\begin{align*}
 \mu_k^y(\delta)(j)&=
 \begin{cases}
 -\delta(k),&j=k,\\
 \delta(j)+[-b_{kj}]_+\delta(k)+b_{kj}[\delta(k)]_+,&j\ne k,
 \end{cases}\\
 \mu_k^{\check y}(\dtc)(j)&=
 \begin{cases}
 -\dtc(k),&j=k,\\
 \dtc(j)+[b_{kj}]_+\dtc(k)-b_{kj}[\dtc(k)]_+,&j\ne k.
 \end{cases}
\end{align*}

\subsection{Canonical decomposition}
The \emph{canonical decomposition}
\[ \delta=\delta_1\oplus\cdots\oplus\delta_s \]
is the decomposition of a general presentation into indecomposable presentation classes. The direct-sum symbol specifies a generic decomposition, not merely equality of weight vectors. For a prescribed sum $\delta=\sum_a\delta_a$, the basic criterion is
\begin{equation} \label{eq:DF} \delta=\bigoplus_a\delta_a\quad\Longleftrightarrow\quad \delta_a\text{ are indecomposable and } \ee(\delta_a,\delta_b)=0\ (a\ne b). \end{equation}
This is \cite[Theorem 4.4]{DF}. Canonical decomposition and multiplicities are transported by extended mutation, using generic preservation under $\tau$ as in Section~\ref{ss:extended-prelim}. General Hom dimensions are additive on these decompositions.

We separate the negative canonical part by writing
\begin{equation} \label{eq:negative-part} \delta=\delta_{\nf}\oplus\bigoplus_i\mneg_i(\delta)(-\mathbf e_i), \qquad \mneg_i(\delta)\geq0. \end{equation}
The decoration of a general representation has dimension $\mneg(\delta)$. In particular $\mneg_i(\delta)$ is not the coordinatewise number $[-\delta(i)]_+$. By \eqref{eq:DF},
\[ \mneg_i(\delta)>0\quad\Longrightarrow\quad\dv(\delta)(i)=0, \]
since $\ee(-\mathbf e_i,\delta_{\nf})=\dim\coker(\delta_{\nf})(i)$.

A weight lies in a \emph{reachable cluster cone} if its canonical decomposition is a nonnegative integral sum of members of one reachable cluster. This is the class of weights of reachable cluster monomials. A weight is called \emph{real} if a general presentation is rigid. Its rigid presentation orbit is open and hence dense, since its normal space is $\E(d,d)$; see \cite[Section 2]{DF}.

\subsection{Extended mutations and general presentations}\label{ss:extended-prelim}
An \emph{extended mutation sequence} is a finite composition of ordinary QP mutations and the decorated Auslander--Reiten operations $\tau^{\pm1}$. We use $w_y$ for its action on weights. An indecomposable weight is \emph{extended-reachable} if an extended mutation sequence takes it to a negative simple. It is \emph{reachable} if such a sequence uses no $\tau$-operation. We reserve these terms for indecomposable weights; a nonnegative sum in one reachable cluster is instead said to lie in a reachable cluster cone.

General presentation families, including general pairs, are preserved by $\tau^{\pm1}$, by \cite[Theorem 3.11]{FeiRank}. Together with ordinary mutation, these operations therefore preserve canonical decompositions, their multiplicities, and the generic values of $\symE$, by Lemma~\ref{lem:symmetric-invariance}.

The weight and coweight operations satisfy
\begin{equation} \label{eq:AR-weight} \tau_y^{-1}\delta=-\delta-\dv(\delta)B, \qquad \tau_{\check y}\dtc=-\dtc+\dv(\dtc)B. \end{equation}
These are also recorded in \cite[(4.7)]{FeiSyzygy}. In particular, the dimension in the second formula belongs to a coweight, not a projective weight.

A \emph{negative normalization} of an indecomposable weight $\eta$ is a pair $(w,i)$ such that $w_y\eta=-\mathbf e_i$. Write $B^w$ and $J^w$ for the exchange matrix and Jacobian algebra after the ordinary mutations in $w$; the operations $\tau^{\pm1}$ do not change them. The row $b_i^w$ is the $i$-th row of $B^w$. Independence of the chosen normalization will be proved in Section~\ref{sec:vectors}.

\subsection{General ranks}
General ranks for quiver representations were studied by Schofield \cite{Sch92}; see also \cite[Lemma 5.1]{FeiRank} for representations with relations. For fixed $X,Y\in\rep J$ and a finite-dimensional linear subspace $V\subseteq\Hom_J(X,Y)$, define
\[ \grk_V(X,Y)(j)=\max_{f\in V}\rank f(j). \]
Each maximum is attained on a nonempty open subset of $V$, and a single general $f$ attains all of them. For $V=\Hom_J(X,Y)$, write $\grk_J(X,Y)$, or $\grk(X,Y)$ when the algebra is fixed.

\begin{lemma}\label{lem:generic-rank}
Let $\pi:V\twoheadrightarrow V'$ be a surjective linear map, and let $A_1,\ldots,A_s$ be linear evaluations on $V'$. The generic ranks of $A_j\pi$ on $V$ equal those of $A_j$ on $V'$, simultaneously for all $j$. Generic evaluated ranks are additive under finite direct sums.
\end{lemma}
\begin{proof}
A rank condition is given by minors, so the maximal-rank locus is nonempty and open. The inverse image of each such locus under $\pi$ is nonempty and open. Their finite intersection is nonempty because $V$ is irreducible. For direct sums the evaluations are diagonal, and one chooses a point in the intersection of the maximal-rank loci of the summands.
\end{proof}

\begin{lemma}\label{lem:generic-defect}
Let $A$ be a finite-dimensional algebra, let $P_-,P_+$ be fixed projectives, and let $p_0\in V=\Hom_A(P_-,P_+)$. On a representation variety of fixed dimension, the function
\[ N\longmapsto \max_{p\in V}\rank p_N^*-\rank (p_0)_N^* \]
is constructible and has a constant value on a nonempty open subset of each irreducible component.
\end{lemma}
\begin{proof}
We choose a basis $p_1,\ldots,p_s$ of $V$. Evaluation of $\sum_a t_ap_a$ is a matrix polynomial in the representation coordinates and linear in $t$. The rank of a general linear combination is at least $r$ precisely when some coefficient of an $r$-minor, viewed as a polynomial in the parameters $t_a$, is nonzero. This is an open condition on $N$. The fixed evaluation at $p_0$ also has open maximal-rank loci. On each irreducible component these two loci meet, and their difference is constant there. Finite rank stratifications prove constructibility. No semicontinuity of the difference is asserted.
\end{proof}

The same argument permits finitely many generality requirements to be imposed simultaneously. It does not choose a representative that is general for infinitely many unrelated weights.

\subsection{Tropical \texorpdfstring{$F$}{F}-polynomials}\label{ss:tropical-F}
For $N\in\rep J$, we define
\[ f_N(\delta)=\max_{L\subseteq N}\pair{\delta}{\dv L},\qquad \fc_N(\delta)=\max_{N\twoheadrightarrow L}\pair{\delta}{\dv L}. \]
The submodule Newton polytope $\mathsf N(N)$ is the convex hull of the submodule dimension vectors. If $\alpha$ is a vertex of $\mathsf N(N)$, then $\Gr_\alpha(N)$ is a point \cite[Theorem 4.2]{FeiComb}. Hence every vertex coefficient of
\[ F_N(\by)=\sum_{\alpha}\chi(\Gr_\alpha(N))\by^\alpha \]
is $1$. In particular, $F_N$ has constant term $1$, and $\mathsf N(N)$ is its Newton polytope. The \emph{quotient $F$-polynomial} is
\[ \check F_N(\by)=\sum_\beta\chi(\Gr^\beta(N))\by^\beta
 =\by^{\dv N}F_N(\by^{-1}), \]
where $\Gr^\beta(N)$ parametrizes quotients of dimension $\beta$, or equivalently their kernels.

For a weight $\delta$, the symbols $f_\delta$ and $\fc_\delta$ denote the tropical functions of a general ordinary part in $\PC(\delta)$. For a coweight $\dtc$, use $f_\dtc$ and $\fc_\dtc$ for the corresponding generic tropical functions. We have
\[ f_N(\delta)-\fc_N(-\delta)=\pair{\delta}{\dv N} =\hom(\delta,N)-\ee(\delta,N). \]
For every finite-dimensional algebra, \cite[Lemma 3.5]{FeiTrop} gives
\begin{equation} \label{eq:tropical-bound} f_N(\delta)\leq\hom(\delta,N). \end{equation}
For an independently variable weight $\delta$ and coweight $\etc$, the \emph{generic pairing} is the equality
\[ f_\etc(\delta)=\hom(\delta,\etc)=\fc_\delta(\etc). \]
It is known when one argument is reachable, and in further cases described in \cite[Theorem 3.22]{FeiTrop} and \cite[Theorem 7.1]{FeiRank}.

For any finite-dimensional algebra and any module $K$, one has
\[ f_K(\sigma)=\hom(\sigma,K)\qquad\text{if }\ee(\sigma,\sigma)=0, \]
by \cite[Theorem 3.6]{FeiTrop}. Thus this equality holds on the nonnegative integral span of the weights of any rigid presentation.

\section{Cluster characters and denominator vectors}\label{sec:cluster}
\subsection{Seeds and initial-seed mutation}
For a skew-symmetric exchange matrix $B$, a seed is a pair $(B,\bx)$ with $\bx=(x_1,\ldots,x_n)$ algebraically independent. We use coefficient-free exchange relations
\[ x_kx_k'=\prod_jx_j^{[b_{kj}]_+}+\prod_jx_j^{[-b_{kj}]_+}, \qquad x_j'=x_j\quad(j\ne k). \]
The two monomials may be interchanged by transposing the skew-symmetric matrix. Matrix mutation is
\[
 b'_{ij}=\begin{cases}
 -b_{ij},&i=k\text{ or }j=k,\\
 b_{ij}+[b_{ik}]_+[b_{kj}]_+-[-b_{ik}]_+[-b_{kj}]_+,&i,j\ne k.
 \end{cases}
\]
The cluster algebra is generated by all cluster variables \cite{FZ1}; the upper cluster algebra is the intersection of the Laurent rings of all reachable seeds \cite{BFZ}.

We write $t$ for an ordinary cluster seed. Mutating a decorated representation changes the initial seed in its character formula. Mutating a cluster variable while keeping its Laurent coordinates fixed is final-seed mutation. We distinguish these operations when using denominator matrices.

\subsection{Generic characters and monomial specialization}
For a general decorated representation $\N=(N,N^-)$ of coweight $\dtc$, we use the generic character
\[ \CC(\dtc)=\bx^{-\dtc}F_N(\by),\qquad y_i=\bx^{b_i}. \]
Equivalently, if $\delta$ labels the same principal component, it is
\[ \check C_{\gen}(\delta)=\bx^{-\delta}\check F_N(\bx^{-b_1},\ldots,\bx^{-b_n}). \]
The character is constant on a suitable dense open subset and respects canonical decomposition. Its mutation compatibility follows from the representation mutation formulas of \cite[Section 5]{DWZ2} and their generic form in \cite[Corollaries 1.11--1.12]{GLFS}. Reachable indecomposables correspond to cluster variables, and reachable cluster cones to cluster monomials \cite[Theorem 5.1 and Lemma 9.2]{DWZ2}. If $\eta$ is a reachable indecomposable weight with corresponding coweight $\check\eta$, we write $X(\eta)=\CC(\check\eta)$ for its cluster variable.

Laurent denominators are used only for nonzero characters; this is automatic for cluster monomials. For a nonzero Laurent polynomial $z$, write
\[ z=\frac{P(\bx)}{\bx^{\den_t(z)}}, \]
where $P$ is not divisible by any $x_i$. Thus $\den_t(x_i)=-\mathbf e_i$. When $\CC(\dtc)\ne0$ and $\delta$ is its corresponding weight, we abbreviate $\den_t(\delta)=\den_t(\CC(\dtc))$. For a generic character we also use the \emph{tropical denominator}
\begin{equation} \label{eq:tropical-denominator} \dT(\delta)(i)=f_N(b_i)+\dtc(i). \end{equation}
The notation $\delta$ in this expression always refers to the weight corresponding to $\dtc$.

\begin{lemma}\label{lem:Newton-specialization} We have the inequality:
\[ \den_t\bigl(\CC(\dtc)\bigr)(i)\leq\dT(\delta)(i). \]
Equality holds whenever the terms of $F_N$ for which $\pair{b_i}{\alpha}$ is maximal do not cancel after the substitution $y_j=\bx^{b_j}$. This is automatic when $B$ has full rank or the coefficients of $F_N$ are nonnegative; in particular it holds for reachable cluster monomials.
\end{lemma}
\begin{proof}
A term $\by^\alpha$ contributes exponent $-\dtc+\alpha B$, whose $i$-th coordinate is $-\dtc(i)-\pair{b_i}{\alpha}$. Thus \eqref{eq:tropical-denominator} gives the smallest possible $x_i$-exponent before collecting equal Laurent monomials. The inequality follows, and equality holds when the terms attaining this value do not cancel. The stated sufficient conditions are immediate; the reachable case follows from positivity \cite{LS}.
\end{proof}

\subsection{Ordered denominator compatibility}
For two cluster variables $x,z$, let $d(x,z)$ be the exponent of $x$ in the denominator of $z$, expressed in a cluster containing $x$. This ordered compatibility degree is precisely the one formalized in \cite[Definition~7.1 and Remark~7.2]{CL}; its independence of the other members of the cluster and its sign alternatives are \cite[Theorem~6.3]{CL}. It gives
\begin{equation} \label{eq:d-sign} d(x,z)=-1\Longleftrightarrow x=z,\qquad d(x,z)=0\Longleftrightarrow x\ne z\text{ and }x,z\text{ lie in a cluster}; \end{equation}
all other values are positive integers. Symmetry of the zero relation does not mean symmetry of the positive values.

For a cluster variable $z$, its \emph{compatibility vector} in a seed $t$ is
\[ \gamma_t^d(z)=\bigl(d(z,x_{1;t}),\ldots,d(z,x_{n;t})\bigr). \]
The first argument is the variable whose denominator exponent is measured. The usual denominator vector of $z$ has the arguments reversed.

\begin{lemma}\label{lem:d-row-mutation}
The vector $\gamma_t^d(z)$ satisfies the tropical maximum rule
\[ \gamma_{t'}^d(z)(k) =\max\left\{\sum_j[b_{kj}]_+\gamma_t^d(z)(j), \sum_j[-b_{kj}]_+\gamma_t^d(z)(j)\right\} -\gamma_t^d(z)(k), \]
with the other coordinates unchanged. Its value is $-\mathbf e_i$ in a seed where $z=x_{i;t}$.
\end{lemma}
\begin{proof}
This is the denominator initialization and final-seed recurrence \cite[(7.6)--(7.7)]{FZ4}, with a seed containing $z$ fixed as the initial seed. Our row convention interchanges the two monomials in their formula and hence does not change the maximum.
\end{proof}

\section{Local Hessian defects}\label{sec:den}
\subsection{The local Hessian defect}
Fix $i$ and define the affine presentation space
\[ V_i=\Hom_J(P_{\inn}^i,P_{\out}^i)=\PHom_J(b_i). \]
The equality uses the absence of oriented $2$-cycles at $i$. Evaluation on a fixed module $N$ is the linear map
\[ V_i\longrightarrow\Hom_\kk(N_{\out}^i,N_{\inn}^i), \qquad p\longmapsto p_N^*. \]
Its image is denoted by $\HH_i(N)$.

\begin{definition}\label{def:local}
The \emph{local Hessian defect} of $N$ at $i$ is
\[ h_i(N)=r_i(N)-\rank\Gamma_i, \qquad\Gamma_i=(p_i)_N^*. \]
Here $r_i(N)=\max_{p\in V_i}\rank p_N^*$ is the maximal rank in the evaluated presentation space $\HH_i(N)$.
\end{definition}

The map $\Gamma_i$ is a member of $\HH_i(N)$, so $h_i(N)\geq0$. In general $\HH_i(N)$ is a proper linear subspace of all maps between the two vector spaces. All ranks in Definition~\ref{def:local} are taken in the specified image space.

\begin{theorem}\label{thm:local}
For every decorated representation $\N=(N,N^-)$,
\begin{equation} \label{eq:local-identity} \hom(b_i,N)+\dtc_\N(i) =\dim N(i)-\dim N^-(i)-h_i(N). \end{equation}

\end{theorem}
\begin{proof}
For general $p\in V_i$, \eqref{eq:presentation-hom} gives
\[ \hom(b_i,N)=\dim N_{\out}^i-r_i(N). \]
On the other hand, \eqref{eq:coindex-local} gives
\[ \dtc_\N(i)=\dim N(i)-\dim N_{\out}^i+\rank\Gamma_i-\dim N^-(i). \]
Adding cancels $\dim N_{\out}^i$ and leaves \eqref{eq:local-identity}.
\end{proof}

\begin{proposition}\label{prop:Hom-defect}
For every finite-dimensional $J$-module $N$, we have
\begin{align}
 h_i(N)&=\hom(\rad P_i,N)-\hom(b_i,N),\label{eq:Hom-defect}\\
 &=\dim N(i)-\hom(S_i,N)+\ext^1(S_i,N)-\hom(b_i,N).
 \label{eq:HomExt-defect}
\end{align}
The function $h_i$ is additive on finite direct sums and has a general value on every irreducible representation component.
\end{proposition}
\begin{proof}
The special presentation has cokernel $\rad P_i$, by \eqref{eq:simple-resolution}. Hence its evaluated rank is
\[ \dim N_{\out}^i-\hom(\rad P_i,N), \]
whereas the general evaluated rank is $\dim N_{\out}^i-\hom(b_i,N)$. Subtraction gives \eqref{eq:Hom-defect}. Applying $\Hom_J(-,N)$ to
\[ 0\longrightarrow\rad P_i\longrightarrow P_i\longrightarrow S_i\longrightarrow0 \]
and using projectivity of $P_i$ gives \eqref{eq:HomExt-defect}. Generic evaluated rank is additive on a finite direct sum by Lemma~\ref{lem:generic-rank}; the special rank is additive as well. The general value is justified by Lemma~\ref{lem:generic-defect}.
\end{proof}

\begin{corollary}\label{cor:algebra-invariance}
The function $h_i(N)$ depends only on $J$, its vertex idempotents $(e_j)_{j\in Q_0}$, and $N$. It is preserved by algebra isomorphisms that send each $e_j$ to the corresponding vertex idempotent. Any map between the prescribed terms $P_{\inn}^i,P_{\out}^i$ with cokernel $\rad P_i$ has the same evaluated rank as $p_i$, for every $N$.
\end{corollary}
\begin{proof}
The algebra $J$ with its vertex idempotents determines the Gabriel quiver, hence $b_i$, both projective terms, and $\rad P_i$. Formula~\eqref{eq:Hom-defect} is therefore intrinsic. For any presentation with the stated cokernel, the kernel of the evaluated map is $\Hom(\rad P_i,N)$, so its rank is the same.

\end{proof}

Thus $h_i(N)$ is the difference between $\hom(\rad P_i,N)$ and the generic value $\hom(b_i,N)$. No comparison between the two cokernels themselves is needed.

\begin{lemma}\label{lem:two-forms}
If $\delta$ and $\dtc$ label the same principal component, then
\[ q_i(\delta):=\hom(\delta,-b_i)+\delta(i) =\hom(b_i,\dtc)+\dtc(i). \]
For a general $\N=(N,N^-)$ in this component,
\[ q_i(\delta)=\dim N(i)-\mneg_i(\delta)-h_i(N). \]
\end{lemma}
\begin{proof}
For general $p\in\PHom(b_i)$, the Nakayama image $\nu p:I_{\inn}^i\to I_{\out}^i$ is a general injective presentation of coweight $-b_i$. Presentation Auslander--Reiten duality gives
\[ \hom(\delta,-b_i)=\ee(b_i,\delta) =\hom(b_i,\dtc)-\pair{b_i}{\dv(\delta)}. \]
Since $\dtc(i)=\delta(i)-\pair{b_i}{\dv(\delta)}$, the first equality follows. Apply Theorem~\ref{thm:local} and identify the decoration using \eqref{eq:negative-part}.
\end{proof}

For a general ordinary part $N$ of weight $\delta$, write $\boldsymbol h(\delta)=(h_i(N))_{i\in Q_0}$.

\subsection{A simultaneous choice of potential}\label{ss:very-general-potential}
We now choose one potential for which the zero-support rank condition holds throughout its mutation class.

\begin{theorem}\label{thm:simultaneous-zero-support}
Suppose $Q$ admits a Jacobi-finite potential. There is a countable family $\mathcal F_Q$ of nonzero polynomials in finitely many potential coefficients such that
\[ \mathcal V_Q=\bigcap_{F\in\mathcal F_Q}\{W:F(W)\ne0\} \]
is nonempty. Every $W\in\mathcal V_Q$ is Jacobi-finite and nondegenerate. For every mutation-equivalent QP and every dimension vector $\alpha$ with $\alpha(i)=0$, the general representation $N$ in each irreducible component satisfies $h_i(N)=0$.
\end{theorem}

Let $\mathcal P(Q)$ be the space of potentials modulo cyclic equivalence. A \emph{finite-jet open set} is the inverse image of a Zariski-open set under truncation to finitely many potential coefficients. We use the complete path-algebra and closed Jacobian-ideal conventions of \cite[Sections~2--3]{DWZ1}; write $\mathfrak m$ for the arrow ideal and $\mathfrak J(W)$ for the closed Jacobian ideal. We also use the elementary consequence that if $W_0$ is Jacobi-finite, then some nonempty principal finite-jet neighbourhood of $W_0$ consists of Jacobi-finite potentials: choose $N$ with $\mathfrak m^N\subseteq\mathfrak J(W_0)$, impose the corresponding maximal-rank condition modulo $\mathfrak m^{N+1}$, and use closedness to lift the inclusion.

\begin{lemma}\label{lem:zero-support-open}
Suppose $Q$ admits a Jacobi-finite potential. For every vertex $i$ and dimension vector $\alpha$ with $\alpha(i)=0$, there is a nonempty finite-jet open set $U_{i,\alpha}(Q)$ of Jacobi-finite potentials such that the general representation $N$ in each irreducible component of $\rep_\alpha J(Q,W)$ satisfies $h_i(N)=0$.
\end{lemma}
\begin{proof}
For $\alpha=0$ the rank assertion is immediate, so take any nonempty finite-jet open set of Jacobi-finite potentials as above. Set $d=\sum_j\alpha(j)>0$. Work first in the finite-type family defined by the Jacobian equations and by the vanishing of all paths of length $d$. Every nilpotent representation of total dimension $d$ satisfies the latter equations. Its relations and Hessian evaluation depend only on potential coefficients through degree $d+1$.

Decompose $W=W^{(0)}+W^{(1)}+W^{(\geq2)}$ according to the number of visits to $i$. A cycle visiting $i$ once has the unique form $abp$ up to cyclic permutation, where $a:j\to i$, $b:i\to k$, and $p:k\to j$ avoids $i$. Since $\alpha(i)=0$, the representation equations depend only on $W^{(0)}$, and
$\Gamma_i(N)=\sum_{a,b,p}c_{abp}A_{abp}(N)$.
Here $A_{abp}(N)$ is the path operator $N(p)$ in the indicated matrix entry; only paths of length less than $d$ are needed. Terms visiting $i$ at least twice act by zero after the second derivative. The matrices $A_{abp}(N)$ span the evaluated presentation space: paths through $i$ act by zero, and every remaining path entry occurs in this list.

Let $K_0$ be the rational function field in the coefficients of $W^{(0)}$. Over $\overline K_0$, choose an irreducible component $Z$ of the nilpotent representation variety on the prescribed support, with function field $L$. This choice is made before adjoining the coefficients $c_{abp}$. At the generic point of $Z$, let $D_Z(t)\in L[t_{abp}]$ be a nonzero maximal minor of $\sum t_{abp}A_{abp}(N)$; when the maximal rank is zero, take $D_Z=1$. Adjoin the $c_{abp}$ as algebraically independent elements over $L$. The substitution
$L[t_{abp}]\longrightarrow L(c_{abp}),\ t_{abp}\longmapsto c_{abp}$,
is injective, so $D_Z(c)\ne0$. The coefficients of $W^{(\geq2)}$ affect neither the representation equations nor the evaluated Hessian, and adjoining them preserves this nonvanishing.

The rank-deficient locus is constructible by the minor argument of Lemma~\ref{lem:generic-defect}. The components over $\overline K_0$ remain geometrically integral after the coefficient extensions. Thus the rank-deficient locus is nowhere dense in every component of the geometric generic fibre. Its closure in the total finite-type family is still nowhere dense on that fibre. By \cite[Lemma 37.24.3, Tag 054X]{Stacks}, after restricting the coefficient space to a nonempty open subset, this closure is nowhere dense in every fibre. Intersect with a nonempty finite-jet open set of Jacobi-finite potentials. Since the finite-jet coefficient space is irreducible, the intersection is nonempty. On it every representation is nilpotent and the evaluated path space is exactly $\HH_i(N)$, proving the claim.
\end{proof}

\begin{proof}[Proof of Theorem~\ref{thm:simultaneous-zero-support}]
Choose a nonempty principal finite-jet open set of Jacobi-finite potentials as above. By \cite[Corollary~7.4 and its proof]{DWZ1}, nondegeneracy is likewise imposed by countably many nonzero polynomial conditions; we shall use the same countability argument over $\mathbb C$. Jacobi-finiteness is preserved by mutation by \cite[Corollary~6.6]{DWZ1}.

For each mutation sequence with endpoint $Q'$, each vertex $i$, and each dimension vector $\alpha$ with $\alpha(i)=0$, take a nonempty principal finite-jet open subset of $U_{i,\alpha}(Q')$ from Lemma~\ref{lem:zero-support-open}. The regular mutation construction in the proof of \cite[Proposition~7.3]{DWZ1}, applied backwards along the sequence, pulls this condition back to a nonempty principal finite-jet open subset of $\mathcal P(Q)$. On the Jacobi-finite locus the condition is intrinsic by Corollary~\ref{cor:algebra-invariance}.

There are only countably many triples consisting of a mutation sequence, a vertex, and a dimension vector. Add the corresponding pullback polynomials to the Jacobi-finiteness condition and to the DWZ nondegeneracy family. The countability argument of \cite[Corollary~7.4]{DWZ1} gives a common nonvanishing point, hence the required set $\mathcal V_Q$.
\end{proof}

From now on we fix $W\in\mathcal V_Q$. By Theorem~\ref{thm:simultaneous-zero-support}, the same choice works in every mutation-equivalent QP. We do not repeat this standing hypothesis; statements valid for every Jacobi-finite nondegenerate potential say so explicitly.

\subsection{The local radical bound}
We next bound the defect when the general representation is nonzero at the chosen vertex.

\begin{lemma}\label{lem:source-monotone}
For a finite-dimensional basic algebra $A$, a module $K$, a weight $\delta$, and $r\in\Z_{\geq0}$,
\[ \hom_A(\delta-r \mathbf e_i,K)\leq\hom_A(\delta,K). \]
The inequality also holds for a general target in a fixed presentation class. Consequently, over a Jacobian algebra,
\begin{equation} \label{eq:q-monotone} q_i(\delta)\geq q_i(\delta-r \mathbf e_i)+r. \end{equation}
\end{lemma}
\begin{proof}
We put $P_-=P([-\delta]_+)$, $P_+=P([\delta]_+)$, and enlarge the presentation to
\[ (d,u):P_-\oplus P_i^{\oplus r}\longrightarrow P_+. \]
For all $d,u,K$, the space $\Hom((d,u),K)$ is the subspace of $\Hom(d,K)$ defined by the additional equation $fu=0$. Taking generic values and applying Lemma~\ref{lem:cancellation} proves the first assertion. The same reasoning in the product with a target parameter space proves the second. Apply it to the coweight $-b_i$ and use the definition of $q_i$ to obtain \eqref{eq:q-monotone}.
\end{proof}

Write
\[ A_i=e_iJe_i,\qquad \mathfrak r_i=\rad A_i,\qquad \kappa_i(N)=\dim_\kk\bigl(N(i)/N(i)\mathfrak r_i\bigr). \]
Since $A_i$ is local with residue field $\kk$, Nakayama's lemma identifies $\kappa_i(N)$ with the minimal number of generators of $N(i)$ as a right $A_i$-module. This is not Schur reduction by $\rad\End_J(N)$.

\begin{lemma}\label{lem:first-vanishing}
Let $N_\delta$ be a general ordinary cokernel, and suppose $N_\delta(i)\ne0$. Then
\[ s=\min\{r\in\Z_{\geq0}:\dv(\delta-r \mathbf e_i)(i)=0\} \]
exists, equals the general value $\kappa_i(N_\delta)$, and satisfies
$\mneg_i(\delta-s \mathbf e_i)=0$.
\end{lemma}
\begin{proof}
For a general $d:P_-\to P_+$, put $N=\coker d$. Projectivity makes
\[ \Hom(P_i^{\oplus r},P_+)\longrightarrow\Hom(P_i^{\oplus r},N) \]
surjective. Thus the extra general columns induce a general map $g:P_i^{\oplus r}\to N$, with $\coker(d,u)=\coker g$. At $i$ this cokernel is
\[ N(i)/(v_1A_i+\cdots+v_rA_i), \]
where the $v_a$ are general elements of $N(i)$. It vanishes exactly when their residue classes span $N(i)/N(i)\mathfrak r_i$, or $r\geq\kappa_i(N)$.

The quantity $\kappa_i(N)$ has a general value: it is the dimension of
\[ P_+(i)/(\im d(i)+P_+(i)\mathfrak r_i), \]
with the second subspace fixed. The enlarged parameter space is a product, so the generic choices of $d$ and the extra columns can be made together. Generic cancellation proves the asserted formula for $s$, and Nakayama gives $s\geq1$.

Set $\sigma=\delta-s \mathbf e_i$. If $-\mathbf e_i$ occurred in its canonical decomposition, removing one copy would give a general presentation of weight $\sigma+\mathbf e_i$ with the same ordinary cokernel dimension. This would force $\dv(\delta-(s-1)\mathbf e_i)(i)=0$, contrary to minimality. Hence the negative multiplicity at the first vanishing weight is zero.
\end{proof}

\begin{theorem}\label{thm:local-radical-bound}
For every weight $\delta$ and its general ordinary part $N$, we have
\[ 0\leq h_i(N)\leq\dim_\kk\bigl(N(i)\mathfrak r_i\bigr) =\dim N(i)-\kappa_i(N). \]
Equivalently,
\[ \hom(b_i,N)+\hom(S_i,N) \geq\ext^1(S_i,N)+\kappa_i(N). \]
\end{theorem}
\begin{proof}
The ordinary principal components are irreducible components, by Section~\ref{ss:principal-components}, so Theorem~\ref{thm:simultaneous-zero-support} applies to their general members. If $N(i)=0$, it gives $h_i(N)=0$. Otherwise take $s$ from Lemma~\ref{lem:first-vanishing} and put $\sigma=\delta-s\mathbf e_i$. Then $\dv(\sigma)(i)=\mneg_i(\sigma)=h_i(N_\sigma)=0$, so Lemma~\ref{lem:two-forms} gives $q_i(\sigma)=0$. Hence
\[ q_i(\delta)\geq s=\kappa_i(N). \]
Canonical orthogonality gives $\mneg_i(\delta)=0$ because $N(i)\ne0$. Now Lemma~\ref{lem:two-forms} yields $\dim N(i)-h_i(N)=q_i(\delta)\geq\kappa_i(N)$. Formula~\eqref{eq:HomExt-defect} gives the equivalent inequality.
\end{proof}

\begin{corollary}\label{cor:local-signs}
For every weight $\delta$, we have the following
\[
\begin{aligned}
 q_i(\delta)>0&\Longleftrightarrow\dv(\delta)(i)>0,\\
 q_i(\delta)<0&\Longleftrightarrow\mneg_i(\delta)>0,\\
 q_i(\delta)=0&\Longleftrightarrow\dv(\delta)(i)=\mneg_i(\delta)=0.
 \end{aligned}
\]
\end{corollary}
\begin{proof}
If $\dv(\delta)(i)>0$, Theorem~\ref{thm:local-radical-bound} gives $q_i(\delta)\geq\kappa_i(N)\geq1$, and canonical orthogonality gives $\mneg_i(\delta)=0$. If $\dv(\delta)(i)=0$, Theorem~\ref{thm:simultaneous-zero-support} gives $h_i(N)=0$, so $q_i(\delta)=-\mneg_i(\delta)$.
\end{proof}

\subsection{Cycle actions and vanishing of the defect}\label{ss:cycle-actions}
We put $A_i=e_iJe_i$ and $\mathfrak r_i=\rad A_i$. For $c\in\mathfrak r_i$, let
\[ R_c:N(i)\longrightarrow N(i),\qquad v\longmapsto vc. \]
The following elementary rank estimate compares these cycle actions with the Hessian defect.

\begin{lemma}\label{lem:rank-compression}
Let $\mathscr L\subseteq\Hom_\kk(X,Y)$ be a linear subspace and $H\in\mathscr L$. If $a:U\to X$ and $b:Y\to Z$ satisfy $Ha=bH=0$, then
\[ \rank(bGa)\leq\max_{F\in\mathscr L}\rank F-\rank H \qquad(G\in\mathscr L). \]
\end{lemma}
\begin{proof}
The map $bGa$ factors through the induced map $\bar G:\ker H\to\coker H$. Put $r=\rank H$ and $s=\rank\bar G$. In bases with $H=\op{diag}(I_r,0)$, a nonzero $s$-minor of $\bar G$ is the coefficient of $t^s$ in a suitable $(r+s)$-minor of $H+tG$. Hence
\[ r+s\leq\rank_{\kk(t)}(H+tG)\leq\max_{F\in\mathscr L}\rank F. \]
The last inequality follows because every larger minor vanishes identically on $\mathscr L$.
\end{proof}

A related rank degeneration appears in the proof of \cite[Lemma 8.2]{DWZ2}.

\begin{theorem}\label{thm:cycle-action}
For every finite-dimensional $J$-module $N$, we have
\[ \max_{c\in\mathfrak r_i}\rank R_c\leq h_i(N). \]
Thus $h_i(N)=0$ implies $N(i)\mathfrak r_i=0$, without a generality assumption.
\end{theorem}
\begin{proof}
The identities $\alpha_i\Gamma_i=\Gamma_i\beta_i=0$ are \cite[(1.12)]{DWZ2}. Every positive-length cycle at $i$ factors as $bpa$, with $b:i\to k$, $p:k\to j$, and $a:j\to i$. Evaluating the entry $P_j\to P_k$ defined by $p$ between the arrow maps gives $R_{bpa}$. These cycles span $\mathfrak r_i$, so
\[ \{\alpha_iG\beta_i:G\in\HH_i(N)\}=\{R_c:c\in\mathfrak r_i\}. \]
Apply Lemma~\ref{lem:rank-compression} to $H=\Gamma_i$, $a=\beta_i$, $b=\alpha_i$, and $\mathscr L=\HH_i(N)$.
\end{proof}

\begin{corollary}\label{cor:cycle-vanishing}\label{cor:cycle-components}
Let $N$ be the general ordinary part of a principal component. Then
\begin{equation} \label{eq:cycle-vanishing} h_i(N)=0\quad\Longleftrightarrow\quad N(i)\mathfrak r_i=0. \end{equation}
Hence the general local defect at $i$ on an ordinary principal component $Z$ vanishes if and only if every representation in $Z$ is annihilated by $J\mathfrak r_iJ$. For a fixed vertex $i$, the following are equivalent:
\begin{enumerate}
\item the general local defect at $i$ vanishes on every principal component;
\item $h_i(P_i)=0$;
\item $e_iJe_i=\kk e_i$.
\end{enumerate}
In particular, $h_i(N)=0$ when $\dim N(i)=1$.
\end{corollary}
\begin{proof}
The equivalence \eqref{eq:cycle-vanishing} follows from Theorem~\ref{thm:cycle-action} and the local radical bound. Since $N(i)\mathfrak r_i=0$ is a closed condition, it holds for the general member of $Z$ if and only if every member of $Z$ is annihilated by $J\mathfrak r_iJ$.

For the final equivalences, $P_i$ is general of weight $\mathbf e_i$, and $P_i(i)=A_i$ is faithful. Thus $h_i(P_i)=0$ implies $\mathfrak r_i=0$ by Theorem~\ref{thm:cycle-action}; the converse follows from the local radical bound. A nilpotent ideal acts trivially on a one-dimensional module.
\end{proof}

The equivalence \eqref{eq:cycle-vanishing} is a statement about the general member of a principal component; the defect may be larger at a special point.

\begin{corollary}\label{cor:local-Jordan}
Assume $A_i\simeq\kk[u]/(u^\ell)$, where $\ell\geq1$. For the general ordinary part $N$ of a principal component,
\[ h_i(N)=\rank R_u,\qquad \dim N(i)-h_i(N)=\dim_\kk\bigl(N(i)/N(i)u\bigr). \]
The last integer counts the Jordan blocks of $R_u$.
\end{corollary}
\begin{proof}
The lower and upper bounds coincide, since $\dim_\kk\bigl(N(i)\mathfrak r_i\bigr)=\rank R_u$.
\end{proof}

\subsection{Denominator vectors}
\begin{theorem}\label{thm:denominator}
Let $\delta$ label a principal component, and let $\N=(N,N^-)$ be general in it. If $f_N(b_i)=\hom(b_i,N)$, then
\[ \dT(\delta)(i)=q_i(\delta) =\dim N(i)-\mneg_i(\delta)-h_i(N). \]
The same formula holds for the Laurent denominator whenever equality holds in Lemma~\ref{lem:Newton-specialization}. In particular, for every reachable cluster monomial,
\[ \den_t(\delta)=\dv(\delta)-\mneg(\delta)-\boldsymbol h(\delta). \]
Without these equalities one still has
\begin{equation} \label{eq:den-bound} \den_t(\delta)(i)\leq\dT(\delta)(i)\leq q_i(\delta) \leq\dim N(i)-\mneg_i(\delta). \end{equation}
\end{theorem}
\begin{proof}
The equality $f_N(b_i)=\hom(b_i,N)$ turns \eqref{eq:tropical-denominator} into Theorem~\ref{thm:local}. The ordinary denominator assertion is Lemma~\ref{lem:Newton-specialization}. For reachable cluster monomials this equality is known by Section~\ref{ss:tropical-F}, and positivity gives equality in Lemma~\ref{lem:Newton-specialization}. The inequalities follow from that lemma, \eqref{eq:tropical-bound}, and $h_i(N)\geq0$.
\end{proof}
\begin{corollary}\label{cor:regularity}
Let $\N=(N,N^-)$ be general of weight $\delta$, with $\mneg_i(\delta)=0$, and assume $f_N(b_i)=\hom(b_i,N)$. The following are equivalent:
\begin{enumerate}
\item $\dT(\delta)(i)=\dim N(i)$.
\item $h_i(N)=0$.
\item There is a subrepresentation $U\subseteq N$ such that
\begin{equation} \label{eq:DWZ-subrepresentation} \ker\Gamma_i\subseteq U_{\out}^i,\qquad \Gamma_i(U_{\out}^i)=U_{\inn}^i. \end{equation}
\end{enumerate}
If there is no cancellation after the substitution, (1) is also equivalent to $\den_t(\delta)(i)=\dim N(i)$. For an ordinary reachable representation, the equivalence of (1)--(3) with the Laurent-denominator equality is valid for every Jacobi-finite nondegenerate potential.

Under the standing hypothesis, these conditions are further equivalent to $N(i)\mathfrak r_i=0$. If $A_i\simeq\kk[u]/(u^\ell)$, then $\dT(\delta)(i)$ counts the Jordan blocks of $R_u$; so does the corresponding Laurent denominator entry when equality holds in Lemma~\ref{lem:Newton-specialization}.
\end{corollary}
\begin{proof}
Theorem~\ref{thm:denominator} gives $(1)\Leftrightarrow(2)$ and the Laurent denominator assertion. For $(2)\Leftrightarrow(3)$, use the calculation in the proof of \cite[Corollary 5.5]{DWZ2}. For every $U\subseteq N$,
\[ \dim\ker\Gamma_i-\pair{b_i}{\dv U} =\dim\frac{\ker\Gamma_i}{\ker\Gamma_i\cap U_{\out}^i} +\dim\frac{U_{\inn}^i}{\Gamma_i(U_{\out}^i)}. \]
The minimum of the left side over all subrepresentations is
\[ \dim\ker\Gamma_i-f_N(b_i) =\dim\ker\Gamma_i-\hom(b_i,N)=h_i(N). \]
Both terms on the right are nonnegative, and both vanish exactly when \eqref{eq:DWZ-subrepresentation} holds. Corollaries~\ref{cor:cycle-vanishing} and~\ref{cor:local-Jordan} give the remaining assertions.
\end{proof}

\begin{corollary}\label{cor:den-signs}
Let $\N=(N,N^-)$ be general of weight $\delta$. At a vertex $i$ satisfying $f_N(b_i)=\hom(b_i,N)$, we have
\[ \dT(\delta)(i)>0\Longleftrightarrow\dv(\delta)(i)>0, \qquad \dT(\delta)(i)<0\Longleftrightarrow\mneg_i(\delta)>0. \]
If $\dv(\delta)(i)=0$, its value is $-\mneg_i(\delta)$. The same statements hold for Laurent denominators whenever equality holds in Lemma~\ref{lem:Newton-specialization}.
\end{corollary}
\begin{proof}
Combine Theorem~\ref{thm:denominator} with Corollary~\ref{cor:local-signs}.
\end{proof}

\begin{remark}\label{rem:cluster-signs}
For ordinary reachable cluster monomials, the same sign assertions hold without using Theorem~\ref{thm:simultaneous-zero-support}. Their factors have ordered denominator degree $-1$, zero, or positive by \eqref{eq:d-sign}. If a factor has degree zero relative to $x_i$, the two variables occur in a cluster. Their decorated representations have zero symmetric $\E$-invariant by the reachable cluster correspondence \cite{DWZ2,DF}. Since
$\symE(S_i^-,\M)=\dim M(i)$, this gives vanishing at $i$. Conversely, when $M(i)=0$, the bound \eqref{eq:den-bound} and nonnegativity of the denominator exponent give degree zero, unless the factor is $x_i$ itself. Additivity gives the cluster-monomial assertion. 
\end{remark}

\subsection{Frozen variables}
Freeze a set $F\subseteq Q_0$ and allow mutations only outside $F$. We put
\[ \Uca_F=\bigcap_t\kk[x_{j;t}^{\pm1}:j\notin F][x_i:i\in F]. \]
Recall that a weight $\delta$ is \emph{$\mu$-supported} if
$\dv(\delta)(i)=0$ for $i\in F$.

\begin{proposition}\label{prop:frozen}
If $\delta$ is $\mu$-supported, then $\CC(\dtc)\in\Uca_F$. Conversely, if the denominator formula of Theorem~\ref{thm:denominator} holds at the frozen vertices of the initial seed, then
\[ \CC(\dtc)\in\Uca_F \quad\Longleftrightarrow\quad \delta\text{ is $\mu$-supported}. \]
For reachable cluster monomials the denominator formula is automatic.
\end{proposition}
\begin{proof}
Allowed mutations do not change the frozen coordinates of $\dv(\delta)$. The zero character belongs to $\Uca_F$. For a nonzero character whose frozen dimensions vanish, \eqref{eq:den-bound} makes every frozen denominator exponent nonpositive in every allowed seed. This proves the forward assertion. For the converse, the assumed denominator formula concerns a nonzero character; membership gives nonpositive frozen exponents, and Corollary~\ref{cor:den-signs} gives $\dv(\delta)(i)=0$ for every $i\in F$.
\end{proof}

\section{Hessian vectors and extended mutation}\label{sec:vectors}\label{sec:mutation}
\subsection{The numerical pairing}\label{ss:numerical-pairing}
For a row vector $\gamma\in\Z^{Q_0}$, set
\[
\mu_k^x(\gamma)(j)=
 \begin{cases}
 \gamma(j),&j\ne k,\\
 \max\{\pair{[b_k]_+}{\gamma},\pair{[-b_k]_+}{\gamma}\}-\gamma(k),&j=k.
 \end{cases}
\]
This map is involutive with the exchange matrix mutated. Define
\[ \mathcal I_B(\gamma,\delta)=\hom(\delta,\gamma B)-\pair{\gamma}{\delta}, \]
where $\gamma B$ is a coweight. The statements in this subsection hold for every Jacobi-finite nondegenerate potential and every integral weight.

\begin{lemma}\label{lem:intertwining}
Let $B'=\mu_kB$, $\gamma'=\mu_k^x\gamma$, and $\delta'=\mu_k^y\delta$. Then
\[ \gamma'B'=\mu_k^{\check y}(\gamma B), \qquad \mathcal I_{B'}(\gamma',\delta')=\mathcal I_B(\gamma,\delta). \]
\end{lemma}
\begin{proof}
Put
\[ u=\sum_j[b_{kj}]_+\gamma(j),\qquad v=\sum_j[-b_{kj}]_+\gamma(j),\qquad a=u-v. \]
If $\etc=\gamma B$, then $\etc(k)=-a$ and
$\gamma'(k)=v+[a]_+-\gamma(k)$. For $j\ne k$, matrix mutation gives
\begin{align*}
 (\gamma'B')(j)-\etc(j)
 &=-b_{kj}(\gamma(k)+\gamma'(k))
      +[b_{kj}]_+v-[-b_{kj}]_+u\\
 &=-[b_{kj}]_+[a]_++[-b_{kj}]_+[-a]_+\\
 &=[b_{kj}]_+\etc(k)-b_{kj}[\etc(k)]_+.
\end{align*}
At $k$, the value is $-\etc(k)$. This proves the first equality.

Write $t=\delta(k)$. The weight-mutation rule gives
\begin{align*}
 \pair{\gamma'}{\delta'}-\pair{\gamma}{\delta}
 &=-t\bigl(\gamma'(k)+\gamma(k)\bigr)+tv+[t]_+a\\
 &=-t[a]_++[t]_+a\\
 &=[-t]_+[a]_+-[t]_+[-a]_+.
\end{align*}
By \eqref{eq:generic-Hom-mutation} this is also
$\hom(\delta',\gamma'B')-\hom(\delta,\gamma B)$, since $\etc(k)=-a$. Subtracting proves the second equality.
\end{proof}

\begin{remark}
The Fock--Goncharov tropical pairing is
\[ \mathcal I_B^{\mathrm{FG}}(\gamma,\delta) =f_{\gamma B}(\delta)-\pair{\gamma}{\delta}. \]
It is bounded above by $\mathcal I_B(\gamma,\delta)$, with equality whenever $f_{\gamma B}(\delta)=\hom(\delta,\gamma B)$. This pairing is distinct from the $f$-compatibility degree in Remark~\ref{rem:f-compatibility}.
\end{remark}

\subsection{Auslander--Reiten translation}\label{ss:AR-pairing}
We put
\[ \tau_x\gamma=\dv(\gamma B)-\gamma. \]
Here $\gamma B$ is a coweight; the argument $-\gamma B$ in the inverse formula below is a weight.

Presentation Auslander--Reiten duality gives
\begin{equation} \label{eq:AR-Hom-pairing} \mathcal I_B(\tau_x\gamma,\tau_y\delta)=\mathcal I_B(\gamma,\delta), \end{equation}
by \cite[Corollary 4.6, equation (4.9)]{FeiSyzygy}.

\Needspace{12\baselineskip}
\begin{proposition}\label{prop:extended-action}
The map $\tau_x$ is a bijection of $\Z^{Q_0}$, with inverse
\begin{equation} \label{eq:tau-x-inverse} \tau_x^{-1}\gamma=\dv(-\gamma B)-\gamma. \end{equation}
It intertwines with coweight translation through multiplication by $B$, and it commutes with ordinary tropical $X$-mutation when the exchange matrices are changed accordingly. Hence every extended mutation sequence has an invertible action $w_x$, and $\mathcal I_B$ is invariant under extended mutation:
\begin{equation} \label{eq:extended-I} \mathcal I_{B^{w}}(w_x\gamma,w_y\delta)=\mathcal I_B(\gamma,\delta). \end{equation}
\end{proposition}
\begin{proof}
Let $\M=(M,M^-)$ be general of coweight $\dtc_\M=\gamma B$, and write $d=\dv M$ and $\gamma'=d-\gamma$. Its weight is
\[ \delta_\M=\dtc_\M-dB=-\gamma'B. \]
The equality of the projective and injective principal-component dimensions gives $\dv(-\gamma'B)=d$. Thus the right-hand side of \eqref{eq:tau-x-inverse}, evaluated at $\gamma'$, is $\gamma$. Conversely, start with a general presentation of weight $-\gamma B$; its coweight is $(\dv(-\gamma B)-\gamma)B$. The same argument proves the other inverse identity. Formula~\eqref{eq:AR-weight} also gives
\[ (\tau_x\gamma)B=\tau_{\check y}(\gamma B). \]

For commutation with $\mu_k^x$, use Lemma~\ref{lem:intertwining}, \eqref{eq:AR-Hom-pairing}, and the commutation of $\mu_k^y$ with $\tau_y$ from \cite[Proposition 7.10]{DF}. The vectors $\mu_k^x\tau_x\gamma$ and $\tau_x\mu_k^x\gamma$, for the same mutated exchange matrix, have equal pairings with every transformed weight. The pairing separates its first argument, since
$\mathcal I_B(\alpha,-\mathbf e_j)=\alpha(j)$.
For the weight $-\mathbf e_j$, the ordinary part is zero. The transformed weights range over the full lattice, so evaluation at $-\mathbf e_j$ proves equality of the vectors even when $B$ is singular. Equation~\eqref{eq:extended-I} follows from Lemma~\ref{lem:intertwining} and \eqref{eq:AR-Hom-pairing} by composition.
\end{proof}
\subsection{Definition of the Hessian vector}\label{ss:Hessian-definition}
Let $\eta$ be indecomposable and extended-reachable. Choose a negative normalization $(w,i)$, so $w_y\eta=-\mathbf e_i$, and write $J^w$ for the Jacobian algebra after applying the normalization. For each original vertex $j$, set
\[ \N_j=w(S_j^-)=(N_j,N_j^-),\qquad \delta_j=w_y(-\mathbf e_j). \]
The representation $\N_j$ is rigid and is general in the principal component of $\delta_j$. Local defects after the normalization will be denoted by $h_i^{J^w}$.

\begin{definition}\label{def:Hessian-vector}
The \emph{Hessian defect vector} and \emph{Hessian vector} associated to the normalization $w$ are
\[ \Delta^w(\eta)(j)=h_i^{J^w}(N_j),\qquad \gamma^w(\eta)=\dv(\eta)-\mneg(\eta)-\Delta^w(\eta). \]
\end{definition}

Thus the definition uses actual representations obtained from the negative simples and general ranks in one presentation space over $J^w$. A single general presentation in $\PHom_{J^w}(b_i^w)$ attains the required evaluated ranks on all $N_j$ simultaneously. The defect coordinates are nonnegative. The next lemma identifies the dimensions and decorations that occur in the definition.

\begin{lemma}\label{lem:endpoint-dimensions}
For every $j$, we have
\begin{equation} \label{eq:endpoint-dimensions} \dim N_j(i)=\dv(\eta)(j),\qquad \dim N_j^-(i)=m_\eta(-\mathbf e_j)=\mneg_j(\eta). \end{equation}
\end{lemma}
\begin{proof}
Evaluation of a negative presentation gives
\[ \symE(-\mathbf e_i,\delta_j)=\dv(\delta_j)(i),\qquad \symE(\eta,-\mathbf e_j)=\dv(\eta)(j). \]
These numbers agree by simultaneous extended invariance of $\symE$, since
$w_y\eta=-\mathbf e_i$ and $w_y(-\mathbf e_j)=\delta_j$. This proves the first equality. Preservation of canonical multiplicities gives
$\mneg_i(\delta_j)=m_\eta(-\mathbf e_j)$.
Since $\eta$ is indecomposable and extended-reachable, this last multiplicity is one exactly when $\eta=-\mathbf e_j$, and is zero otherwise. The same alternatives describe $\mneg_j(\eta)$. Finally $\N_j$ is general, so its decoration is $\mneg(\delta_j)$.
\end{proof}

\begin{theorem}\label{thm:endpoint-rank}
For every Jacobi-finite nondegenerate QP and every negative normalization $(w,i)$ of $\eta$,
\begin{equation} \label{eq:endpoint-rank} \gamma^w(\eta)(j) =\dim N_j(i)-\dim N_j^-(i)-h_i^{J^w}(N_j) =q_i^{J^w}(\delta_j). \end{equation}
Equivalently,
\begin{equation} \label{eq:gamma-transport} \gamma^w(\eta)=w_x^{-1}(-\mathbf e_i). \end{equation}
\end{theorem}
\begin{proof}
Lemma~\ref{lem:endpoint-dimensions} and the local identity of Lemma~\ref{lem:two-forms} give \eqref{eq:endpoint-rank}. Put $\alpha=w_x^{-1}(-\mathbf e_i)$. The invariant pairing gives
\[ \alpha(j)=\mathcal I_B(\alpha,-\mathbf e_j) =\mathcal I_{B^w}(-\mathbf e_i,\delta_j) =\hom_{J^w}(\delta_j,-b_i^w)+\delta_j(i) =q_i^{J^w}(\delta_j). \]
Thus $\alpha=\gamma^w(\eta)$, proving \eqref{eq:gamma-transport}.
\end{proof}

Both identifications in \eqref{eq:endpoint-dimensions} use only presentation invariants: ordinary dimensions are determined by $\symE$, and decorations by canonical multiplicities.

\subsection{Independence and extended mutation}
\begin{theorem}\label{thm:intrinsic-vector}
For every extended-reachable indecomposable $\eta$, the vectors $\gamma^w(\eta)$ and $\Delta^w(\eta)$ are independent of the negative normalization. Write them as $\gamma_J(\eta)$ and $\Delta_J(\eta)$, or simply $\gamma(\eta)$ and $\Delta(\eta)$. They satisfy
\begin{equation} \label{eq:gamma-def} \gamma(\eta)=\dv(\eta)-\mneg(\eta)-\Delta(\eta), \qquad \Delta(\eta)\geq0, \end{equation}
and
\[ \gamma(-\mathbf e_i)=-\mathbf e_i,\qquad \gamma_{J^{\mu_k}}(\mu_k^y\eta)=\mu_k^x\gamma_J(\eta),\qquad \gamma_J(\tau_y\eta)=\dv(\gamma_J(\eta)B)-\gamma_J(\eta). \]
Here $J^{\mu_k}$ is the Jacobian algebra of $\mu_k(Q,W)$.
\end{theorem}
\begin{proof}
We first take $\eta=-\mathbf e_a$ and any negative normalization $(w,i)$. Lemma~\ref{lem:endpoint-dimensions} gives $N_j(i)=0$ for every $j$, while the negative decorations at $i$ have dimensions $\delta_{aj}$. The QP corresponding to $J^w$ satisfies Theorem~\ref{thm:simultaneous-zero-support}. Each $N_j$ is general in an irreducible component, by Section~\ref{ss:principal-components}, so all $h_i^{J^w}(N_j)$ vanish. Formula~\eqref{eq:endpoint-rank} gives
$\gamma^w(-\mathbf e_a)=-\mathbf e_a$ for every normalization.

We now let $w,v$ be two negative normalizations of $\eta$, with
$w_y\eta=-\mathbf e_i$ and $v_y\eta=-\mathbf e_j$.
Then $u=v\circ w^{-1}$ normalizes $-\mathbf e_i$ to $-\mathbf e_j$. By the preceding paragraph and \eqref{eq:gamma-transport},
$u_x^{-1}(-\mathbf e_j)=-\mathbf e_i$. 
Since $u_x^{-1}=w_xv_x^{-1}$, we obtain
$v_x^{-1}(-\mathbf e_j)=w_x^{-1}(-\mathbf e_i)$,
so Theorem~\ref{thm:endpoint-rank} gives $\gamma^v(\eta)=\gamma^w(\eta)$. The defect is then intrinsic by \eqref{eq:gamma-def}.

If $s$ is one ordinary mutation or one Auslander--Reiten translation, then $w\circ s^{-1}$ is a negative normalization of $s_y\eta$. Hence
\[ \gamma(s_y\eta)=(w\circ s^{-1})_x^{-1}(-\mathbf e_i) =s_xw_x^{-1}(-\mathbf e_i)=s_x\gamma(\eta), \]
which gives the two mutation formulas.
\end{proof}

\begin{corollary}\label{cor:gamma-positive}
Under the same hypothesis, if the indecomposable $\eta$ has nonzero ordinary part $M$, then
\[ 0\leq\gamma(\eta)\leq\dv M,\qquad \gamma(\eta)(j)>0\Longleftrightarrow M(j)\ne0. \]
The map $\eta\mapsto\gamma(\eta)$ is injective on extended-reachable indecomposables.
\end{corollary}
\begin{proof}
An indecomposable decorated representation with nonzero ordinary part has no negative decoration. At a negative normalization, \eqref{eq:endpoint-dimensions} identifies $\dim N_j(i)$ with $\dim M(j)$. Corollary~\ref{cor:local-signs} and \eqref{eq:endpoint-rank} give the support assertion. The upper bound follows from $\Delta\geq0$.

If two vectors agree, we apply a negative normalization for the first weight. Their transformed vectors still agree, and the first is $-\mathbf e_i$. The other weight cannot have nonzero ordinary part, since its vector would be nonnegative. It is a negative simple, necessarily $-\mathbf e_i$. Applying the inverse operation to the weights proves injectivity.
\end{proof}

\begin{remark}
For an arbitrary Jacobi-finite nondegenerate potential, Theorem~\ref{thm:endpoint-rank} remains valid for each chosen negative normalization. In the preceding normalization construction, the comparison between different normalizations uses the very-general zero-support theorem. The radical and sign results in Section~\ref{sec:den} also use that theorem.
\end{remark}

\subsection{Ordinary denominator compatibility}\label{ss:ordinary-denominator}
\begin{theorem}\label{thm:ordinary-vector}
Let $\eta$ be reachable and let $t$ be the current cluster seed. 
For an ordinary negative normalization $w$, we have
\[ \gamma^w(\eta)(j)=d(X(\eta),x_{j;t}). \]
Thus $\gamma^w(\eta)$ and $\Delta^w(\eta)$ are independent of the ordinary normalization for every Jacobi-finite nondegenerate potential.
\end{theorem}
\begin{proof}
At the seed reached by $w$, the character of $w\eta=-\mathbf e_i$ is its $i$-th variable. Under the same initial-seed changes, the character of $w(S_j^-)$ is the original variable $x_{j;t}$. Applying the reachable case of Theorem~\ref{thm:denominator} to $w(S_j^-)$ identifies \eqref{eq:endpoint-rank} with the stated ordered denominator degree. Its independence of the other variables of the seed is \cite[Theorem 6.3]{CL}. Since $\dv(\eta)$ and $\mneg(\eta)$ are fixed, the defect is independent as well.
\end{proof}

We can state the relation with denominator vectors as a matrix identity. For two seeds $t,s$, define
\[ \mathsf D(t,s)_{ij}=d(x_{i;t},x_{j;s}),\qquad \mathsf H(t,s)_{ij}=\gamma_t(x_{j;s})(i), \]
where $\gamma_t(x_{j;s})$ denotes the Hessian vector of the class of $x_{j;s}$, computed in seed $t$. Thus the columns of $\mathsf D(t,s)$ are the denominator vectors of the $s$-variables with respect to $t$, while the columns of $\mathsf H(t,s)$ are their Hessian vectors. Theorem~\ref{thm:ordinary-vector} gives
\[ \mathsf H(t,s)=\mathsf D(s,t)^{\tr}. \]
This identity interchanges the two seeds as well as transposing the matrix. It does not assert $\mathsf D(t,s)=\mathsf D(s,t)^{\tr}$. In our skew-symmetric setting, that additional equality is Reading--Stella's Property D \cite[Section 2]{RS}; when it holds, Hessian and denominator vectors agree. It need not hold in general. For classes that are extended-reachable but not reachable, we use the presentation construction rather than an interpretation by ordinary cluster variables.

The following two cluster variables have the same denominator vector but different Hessian vectors.

\begin{example}\label{ex:equal-denominators}
We consider the quiver with potential:
\[
 \begin{tikzpicture}[baseline=(current bounding box.center),
   >=Stealth, every node/.style={inner sep=2.5pt},
   every path/.style={line width=.45pt}]
  \node (v1) at (0,1.65) {$1$};
  \node (v2) at (2.8,1.65) {$2$};
  \node (v3) at (2.8,0) {$3$};
  \node (v4) at (0,0) {$4$};
  \draw[->] (v1) to[bend left=16] node[above] {$a_1$} (v2);
  \draw[->] (v1) to[bend right=16] node[below] {$a_2$} (v2);
  \draw[->] (v2) -- node[right] {$b$} (v3);
  \draw[->] (v3) to[bend right=16] node[above] {$c_1$} (v4);
  \draw[->] (v3) to[bend left=16] node[below] {$c_2$} (v4);
  \draw[->] (v4) -- node[left] {$d$} (v1);
 \end{tikzpicture}
 \qquad W=a_1bc_1d+a_2bc_2d.
\]

The relations include $bc_1d=bc_2d=da_1b=da_2b=0$. Every path of length four vanishes, so the algebra is finite-dimensional. Every cycle is cyclically equivalent to a path containing one of these relations. The QP is therefore rigid and hence nondegenerate by \cite[Corollary~8.2]{DWZ1}.

Let
\[ \eta=(3,-2,0,1),\qquad\zeta=(0,1,3,-2). \]
The ordinary mutation sequence $(2,4,1,3,4,2,3,4)$, applied from left to right to the current QP, takes $\eta$ to $-\mathbf e_4$. Its image under $1\leftrightarrow3$, $2\leftrightarrow4$ takes $\zeta$ to $-\mathbf e_2$. Direct denominator-matrix mutation gives
\[ \den_t(\eta)=\den_t(\zeta)=(4,6,4,6). \]
The variables are distinct because their reachable weights are distinct. Reverse tropical transport from the respective negative simples gives
\[ \gamma(\eta)=(4,6,4,5),\qquad \gamma(\zeta)=(4,5,4,6). \]
These values follow from the denominator recurrence of Lemma~\ref{lem:d-row-mutation}; the two indicated sequences determine them exactly. No equality of denominator and dimension vectors is used.
\end{example}

\section{Projective and injective realizations}\label{sec:projective}
\subsection{Projective and injective vectors}
Let $d_i\in\PHom_J(b_i)$ be general, and put
\[ C_i=\coker d_i,\qquad K_i=\ker(\nu d_i). \]
The injective presentation $\nu d_i$ has coweight $-b_i$. Every projective and injective is extended-reachable, since $\tau P_i=S_i^-$ and $\tau^{-1}I_i=S_i^-$.

\begin{corollary}\label{cor:projective} We have the following equalities:
\begin{align*}
 \gamma(P_i)&=\mathbf e_i+\dv C_i,\\
 \gamma(I_i)&=\mathbf e_i+\dv K_i,\\
 \Delta(P_i)&=\dv\rad P_i-\dv C_i.
\end{align*}
\end{corollary}
\begin{proof}
The AR normalizations give
\[ \gamma(P_i)=\tau_x^{-1}(-\mathbf e_i)=\mathbf e_i+\dv(b_i),\qquad \gamma(I_i)=\tau_x(-\mathbf e_i)=\mathbf e_i+\dv(-b_i). \]
Subtract the first from $\dv P_i=\mathbf e_i+\dv\rad P_i$. The presentation $p_i$ has cokernel $\rad P_i$, while a general presentation in $\PHom(b_i)$ has cokernel of dimension $\dv C_i$.
\end{proof}

\subsection{Denominator vectors of projectives and injectives}
\begin{proposition}\label{prop:projective-denominator}
\[ q_i(P_j)=\delta_{ij}+\dim K_i(j),\qquad q_i(I_j)=\delta_{ij}+\dim C_i(j). \]
These are respectively $\gamma(I_i)(j)$ and $\gamma(P_i)(j)$. If $f_{P_j}(b_i)=\hom(b_i,P_j)$, respectively $f_{I_j}(b_i)=\hom(b_i,I_j)$, then
\[ \dT(P_j)(i)=\gamma(I_i)(j),\qquad \dT(I_j)(i)=\gamma(P_i)(j). \]
They are Laurent denominator identities when equality holds in Lemma~\ref{lem:Newton-specialization}. In particular, the respective Laurent identity holds when $P_j$, respectively $I_j$, is reachable.
\end{proposition}
\begin{proof}
A projective has weight $\mathbf e_j$, so $q_i(P_j)=\hom(\mathbf e_j,-b_i)+\delta_{ij}=\dim K_i(j)+\delta_{ij}$. An injective has coweight $\mathbf e_j$, giving $q_i(I_j)=\hom(b_i,\mathbf e_j)+\delta_{ij}=\dim C_i(j)+\delta_{ij}$. Apply Corollary~\ref{cor:projective} and Theorem~\ref{thm:denominator} and Lemma~\ref{lem:Newton-specialization}.
\end{proof}

In $q_i(P_j)$ and $q_i(I_j)$, the module stands for its presentation class. The denominator of a projective is compared with an injective Hessian vector, and conversely, with the two vertex indices interchanged. The relation between the two ordered denominator degrees is explained in Section~\ref{ss:ordinary-denominator}.

\subsection{Module realizations by deformation}
The following elementary deformation is needed for the projective quotient and injective submodule statements. It does not require an inclusion between the special and general images.

\begin{proposition}\label{prop:deformation-defect}
Let $A$ be a finite-dimensional algebra, $H:X\to Y$ a module map, and $G$ general in $\Hom_A(X,Y)$. Put $\mathcal O=\kk[[t]]$,
\[ \mathscr Q=\coker(H+tG:X\otimes\mathcal O\to Y\otimes\mathcal O),\qquad \mathscr T=\tors_{\mathcal O}\mathscr Q, \qquad D_H=\mathscr T/t\mathscr T. \]
There is an exact sequence
\begin{equation} \label{eq:deformation-special-fibre} 0\longrightarrow D_H\longrightarrow\coker H \longrightarrow (\mathscr Q/\mathscr T)/t(\mathscr Q/\mathscr T) \longrightarrow0, \end{equation}
and
\[ \dv D_H=\grk_A(X,Y)-\rh H. \]
The last term has the dimension vector of a general cokernel. If $\mathscr K=\ker(H+tG)$, then $\mathscr K/t\mathscr K$ embeds in $\ker H$ with the dimension vector of a general kernel.
\end{proposition}
\begin{proof}
At each vertex, choose a nonzero maximal minor of $G$. It is the leading coefficient of the corresponding minor of $H+tG$, so the rank over $\kk((t))$ is the general rank. It cannot be larger, since all larger minors vanish identically on the Hom-space.

The $\mathcal O$-torsion $\mathscr T$ is an $A$-submodule, and $\mathscr Q/\mathscr T$ is finite and torsion-free over the discrete valuation ring $\mathcal O$, hence free. Tensoring its defining sequence with $\mathcal O/(t)$ gives \eqref{eq:deformation-special-fibre}. At a vertex, Smith normal form has $\grk_A(X,Y)(j)$ nonzero invariant factors. Exactly $\rank H(j)$ are units; each remaining factor contributes one dimension to $\mathscr T(j)/t\mathscr T(j)$. This proves the dimension assertions.

Finally $\mathscr K$ is free and saturated in $X\otimes\mathcal O$: if $tx\in\mathscr K$, torsion-freeness of $Y\otimes\mathcal O$ gives $x\in\mathscr K$. Reduction modulo $t$ therefore embeds it into $\ker H$, with the general kernel dimension.
\end{proof}

\begin{corollary}\label{cor:projective-subquotients}
There are submodules $D_i\subseteq\rad P_i$ and $V_i^-\subseteq I_i$ such that
\[ \dv(P_i/D_i)=\gamma(P_i),\qquad \dv V_i^-=\gamma(I_i). \]
More precisely, there are modules $C_i^0,K_i^0$ of dimensions $\dv C_i,\dv K_i$ and exact sequences
\[ 0\to C_i^0\to P_i/D_i\to S_i\to0,\qquad 0\to S_i\to V_i^-\to K_i^0\to0. \]
\end{corollary}
\begin{proof}
Apply Proposition~\ref{prop:deformation-defect} to $p_i+t d_i$. Its special cokernel is $\rad P_i$, and its general cokernel has dimension $\dv C_i$. Thus $D_i\subseteq\rad P_i$ has dimension $\Delta(P_i)$, and $C_i^0=\rad P_i/D_i$ has dimension $\dv C_i$. Quotienting $0\to\rad P_i\to P_i\to S_i\to0$ gives the first sequence.

Dually, the special kernel of $\nu p_i$ is $I_i/S_i$. The kernel part of the proposition gives $K_i^0\subseteq I_i/S_i$ of dimension $\dv K_i$. Its inverse image in $I_i$ is $V_i^-$, giving the second sequence.
\end{proof}

The modules may depend on the general direction $d_i$; neither canonicity nor isomorphisms $C_i^0\simeq C_i$, $K_i^0\simeq K_i$ are asserted. For a general extended-reachable representation, the coordinate formula of Theorem~\ref{thm:endpoint-rank} does not furnish a single map between modules over the original algebra whose rank defect is $\Delta(\eta)$. No such additional realization is used here.

\section{Compatibility degree and canonical decomposition}\label{sec:compat}

Our standing assumption is that the first weight $\eta$ is indecomposable and extended-reachable; the second weight $\delta$ is arbitrary.

\subsection{Compatibility degree}\label{ss:compatibility}
\begin{definition}\label{def:compatibility}
The \emph{compatibility degree} is
\[ \comp{\eta}{\delta}=\mathcal I_B(\gamma(\eta),\delta) =\hom(\delta,\gamma(\eta)B)-\pair{\gamma(\eta)}{\delta}. \]
\end{definition}

We choose an extended mutation sequence $w$ with $w_y\eta=-\mathbf e_i$, put $\delta'=w_y\delta$, and let $\N=(N,N^-)$ be general of weight $\delta'$ after the normalization. The numerical pairing theorem gives
\begin{equation} \label{eq:compat-boundary} \comp{\eta}{\delta}=q_i(\delta'). \end{equation}

\begin{theorem}\label{thm:pair-identity}
The local defect after the normalization $h_i(N)$ is independent of the negative normalization. Denote it by $h(\eta,\delta)$. Then\begin{equation} \label{eq:pair-defect} \comp{\eta}{\delta}=\symE(\eta,\delta)-m_\eta(\delta)-h(\eta,\delta). \end{equation}
The degree and the defect are additive on canonical decompositions in $\delta$ and invariant under simultaneous extended mutation. The defect is nonnegative, and
\[ \Delta(\eta)(j)=h(\eta,-\mathbf e_j). \]
If $e=\symE(\eta,\delta)>0$, then $h(\eta,\delta)\leq e-1$; if $e=0$, the defect is zero.
\end{theorem}
\begin{proof}
After the normalization, we put $n=\dim N(i)$ and $m=\mneg_i(\delta')$. Evaluation of the negative presentation gives
$\ee(-\mathbf e_i,\delta')=n$ and $\ee(\delta',-\mathbf e_i)=0$.
Extended invariance of $\symE$ and canonical multiplicities identifies
$n=\symE(\eta,\delta)$ and $m=m_\eta(\delta)$. The local formula is $q_i(\delta')=n-m-h_i(N)$, proving \eqref{eq:pair-defect}. Its other terms are intrinsic, so $h_i(N)$ is independent of the normalization.

The bound and vanishing assertion follow from Theorem~\ref{thm:local-radical-bound}. Generic Hom dimensions are additive on canonical decompositions. The special evaluated rank is additive, and one general presentation in $\PHom(b_i)$ attains maximal rank on every member of a finite direct sum; hence the defect is additive as well. Transporting both weights and the chosen normalization leaves the calculation at the negative normalization unchanged. Finally take $\delta=-\mathbf e_j$, and let $\M=(M,M^-)$ be the general decorated representation of weight $\eta$. The degree is $\gamma(\eta)(j)$, its symmetric invariant is $\dim M(j)$, and its multiplicity is $\dim M^-(j)$. Equation~\eqref{eq:gamma-def} proves the last identification.
\end{proof}

\begin{proposition}\label{prop:compatibility}
For reachable indecomposables $\eta,\delta$, we have
\[ \comp{\eta}{\delta}=d(X(\eta),X(\delta)). \]
\end{proposition}
\begin{proof}
We choose an ordinary negative normalization of $\eta$. At its final seed, $X(\eta)$ is the initial variable $x_i$. The reachable case of Theorem~\ref{thm:denominator} identifies $q_i(w_y\delta)$ with the exponent of $x_i$ in the denominator of $X(\delta)$. Now use \eqref{eq:compat-boundary}.
\end{proof}

\begin{corollary}\label{cor:symmetric-compatibility}
Let $\eta\ne\delta$ be extended-reachable indecomposable presentation classes. Then
\[ \comp{\eta}{\delta}=0 \quad\Longleftrightarrow\quad \symE(\eta,\delta)=0 \quad\Longleftrightarrow\quad \comp{\delta}{\eta}=0. \]
These conditions are equivalent to rigidity of the direct sum of their presentations. When both classes are reachable, they are also equivalent to $X(\eta)$ and $X(\delta)$ belonging to a common cluster.
If $\symE(\eta,\delta)=1$, then
\[ \comp{\eta}{\delta}=\comp{\delta}{\eta}=1, \qquad h(\eta,\delta)=h(\delta,\eta)=0. \]
In particular, both ordered degrees of an ordinary exchange pair are one.
\end{corollary}
\begin{proof}
Distinct indecomposable classes have zero multiplicity in one another. Thus Theorem~\ref{thm:pair-identity} gives
$\comp{\eta}{\delta}=\symE(\eta,\delta)-h(\eta,\delta)$.
If the symmetric invariant is zero, then the defect is zero. If it is positive, then $h(\eta,\delta)\leq\symE(\eta,\delta)-1$, so the degree is positive. The same argument with the two classes interchanged proves the equivalences. Since both presentations are rigid, vanishing of both $\ee(\eta,\delta)$ and $\ee(\delta,\eta)$ is equivalent to rigidity of their direct sum. In the reachable case, apply Proposition~\ref{prop:compatibility} and \cite[Theorem 6.3]{CL}.

If the symmetric invariant is one, the defect is zero and \eqref{eq:pair-defect} gives both degrees one. For an ordinary exchange pair, choose a seed in which the pair is exchanged. The corresponding presentations are the negative simple at the exchanged vertex and the ordinary simple at that vertex; their symmetric invariant is one. Simultaneous mutation preserves it.
\end{proof}

\begin{remark}[Comparison with $f$-compatibility]\label{rem:f-compatibility}
For ordinary cluster variables $X(\eta),X(\delta)$, let
$c_f(\eta,\delta)=(X(\eta)\parallel X(\delta))_f$ be the \emph{$f$-compatibility degree} of \cite[Definition 4.9]{FuGyoda}. Then
\[ c_f(\eta,\delta)=\symE(\eta,\delta) =[\comp{\eta}{\delta}]_++h(\eta,\delta). \]
After an ordinary normalization $w_y\eta=-\mathbf e_i$, with
$w(\mathcal M_\delta)=(N,N^-)$, the first equality is the standard identity
$c_f(\eta,\delta)=\dim N(i)=\symE(\eta,\delta)$. If $\eta\ne\delta$, then $m_\eta(\delta)=0$, and Theorem~\ref{thm:pair-identity} gives
$c_f(\eta,\delta)=\comp{\eta}{\delta}+h(\eta,\delta)$.
The degree is nonnegative by the same theorem, while for $\eta=\delta$ the formula is immediate from $\comp{\eta}{\eta}=-1$ and rigidity. Hence
$c_f(\eta,\delta)=\symE(\eta,\delta)=[\comp{\eta}{\delta}]_++h(\eta,\delta)$.
For $\eta\ne\delta$, Proposition~\ref{prop:compatibility} therefore gives
\[ h(\eta,\delta)=c_f(\eta,\delta)-d(X(\eta),X(\delta)), \]
and in the same normalization Corollary~\ref{cor:cycle-vanishing} gives
\[ c_f(\eta,\delta)=[\comp{\eta}{\delta}]_+ \quad\Longleftrightarrow\quad N(i)\rad(e_iJ^we_i)=0. \]
\end{remark}

\begin{remark}[Exchangeability and asymmetry]\label{rem:exchange-counterexamples}
The converse of the exchange-pair assertion in Corollary~\ref{cor:symmetric-compatibility} does not hold for denominator degrees. Fu--Gyoda \cite[Examples 6.13 and 6.14]{FuGyoda} give nonexchangeable pairs $x,z$ with
\[ d(x,z)=d(z,x)=1, \qquad (x\parallel z)_f=(z\parallel x)_f=2. \]
Under the identification in Remark~\ref{rem:f-compatibility}, each ordered Hessian defect is one. Their Example 4.21 gives a pair with
\[ d(x,z)=2,\qquad d(z,x)=1. \]
Thus degree one in one order does not force degree one in the other. These examples concern the denominator degree, not the $f$-compatibility degree.
\end{remark}

\subsection{Canonical decomposition}\label{ss:canonical-decomposition}
\begin{theorem}\label{thm:multiplicity}
For every extended-reachable indecomposable $\eta$ and every weight $\delta$,
\[ m_\eta(\delta)=[-\comp{\eta}{\delta}]_+. \]
Moreover,
\[
\begin{aligned}
 \comp{\eta}{\delta}\leq0&\Longleftrightarrow\symE(\eta,\delta)=0,\\
 \comp{\eta}{\delta}\geq0&\Longleftrightarrow m_\eta(\delta)=0.
 \end{aligned}
\]
In particular, positive degree is equivalent to a positive symmetric $\E$-invariant, and negative degree is equivalent to occurrence of $\eta$ as a canonical summand.
Equivalently, the defect identity can be written
\begin{equation} \label{eq:positive-compatibility} \symE(\eta,\delta)=[\comp{\eta}{\delta}]_++h(\eta,\delta). \end{equation}
\end{theorem}
\begin{proof}
Use \eqref{eq:compat-boundary} after a negative normalization $(w,i)$ of $\eta$. Put $\delta'=w_y\delta$, choose a general decorated representation $\N=(N,N^-)$ of weight $\delta$, and write $n=\dim N(i)$ and $m=\mneg_i(\delta)$. If $n>0$, canonical orthogonality gives $m=0$, and Corollary~\ref{cor:local-signs} gives $q_i(\delta')>0$. If $n=0$, the zero-support theorem gives $h_i(N)=0$, so $q_i(\delta')=-m$. Since $n=\symE(\eta,\delta)$ and $m=m_\eta(\delta)$, these two cases prove the multiplicity and sign assertions. Finally, substitute the multiplicity formula into \eqref{eq:pair-defect} and use $c+[-c]_+=[c]_+$ to obtain \eqref{eq:positive-compatibility}.
\end{proof}

The theorem applies to arbitrary presentation weights, including those with imaginary canonical summands.

\begin{corollary}\label{cor:compatible-zero}
If $\symE(\eta,\delta)=0$, then $h(\eta,\delta)=0$ and $\comp{\eta}{\delta}=-m_\eta(\delta)$. In particular this holds whenever $\eta$ occurs in the canonical decomposition of $\delta$.
\end{corollary}
\begin{proof}
Use Theorem~\ref{thm:pair-identity}. If $\eta$ occurs, rigidity and canonical orthogonality give $\symE(\eta,\delta)=0$.
\end{proof}

For acyclic $Q$, let $E_Q=(\delta_{ij}-\#\{i\to j\})_{i,j\in Q_0}$ be the Euler matrix of $Q$, so $B=E_Q^{\tr}-E_Q$.
For dimension vectors $\alpha,\beta$, let $m^Q_\beta(\alpha)$ denote the multiplicity of $\beta$ in Kac's canonical decomposition of $\alpha$ \cite{Kac}. In the following formula, the arguments of $\hom$ are a presentation weight and an injective coweight, not two dimension vectors.

\begin{corollary}\label{cor:acyclic-multiplicity}
Let $\beta$ be a real Schur root of an acyclic quiver $Q$, and let $\alpha\in\Z_{\geq0}^{Q_0}$. Then
\begin{equation} \label{eq:acyclic-multiplicity} m^Q_\beta(\alpha) =\left[\alpha E_Q\beta^{\tr} -\hom_{\kk Q}(\alpha E_Q,\beta B)\right]_+. \end{equation}
Here $\beta B$ is regarded as a coweight.
\end{corollary}
\begin{proof}
Over the hereditary algebra $\kk Q$, a representation of dimension $\alpha$ has a projective presentation of weight $\alpha E_Q$. The correspondence between general representations and general presentations preserves canonical decomposition \cite[Section 2]{DF}; thus
\[ m^Q_\beta(\alpha)=m_{\beta E_Q}(\alpha E_Q). \]
A real Schur root is represented by a reachable Schur representation $M_\beta$, by the acyclic cluster-variable correspondence \cite[Theorem 4]{CK}. Put $\eta=\beta E_Q$. Since $\End(M_\beta)=\kk$, Theorem~\ref{thm:Schur-comparison} gives $\gamma(\eta)=\beta$. Therefore
\[ \comp{\eta}{\alpha E_Q} =\hom_{\kk Q}(\alpha E_Q,\beta B) -\alpha E_Q\beta^{\tr}. \]
Theorem~\ref{thm:multiplicity} proves \eqref{eq:acyclic-multiplicity}. There is no additional potential hypothesis in this case: an acyclic quiver has only the zero potential, and the nonempty set $\mathcal V_Q$ is therefore $\{0\}$.
\end{proof}

The formula computes the multiplicity of a specified real Schur root, even when the canonical decomposition of $\alpha$ also contains imaginary roots. It does not by itself determine the complete decomposition; see Problem~\ref{prob:canonical-algorithm}.

\subsection{Schur reductions and denominator vectors}
Let $M$ be an ordinary indecomposable representation. Its \emph{positive} and \emph{negative Schur reductions} are
\[ \rho_+(M)=M\Big/\sum_{f\in\rad\End_J(M)}\im f, \qquad \rho_-(M)=\bigcap_{f\in\rad\End_J(M)}\ker f; \]
see \cite[Definition 4.9]{FeiProjection}.

\begin{theorem}\label{thm:Schur-comparison}
For every ordinary extended-reachable indecomposable $M$, we have
\begin{equation} \label{eq:Schur-comparison} \dv\rho_+(M)\leq\gamma(M),\qquad \dv\rho_-(M)\leq\gamma(M)\leq\dv M. \end{equation}
In particular, if $\End_J(M)=\kk$, then $\gamma(M)=\dv M$ and $\Delta(M)=0$.
\end{theorem}
\begin{proof}
We first make an observation about a rigid completion. Let $E=\bigoplus_{b=1}^nE_b$ be a basic maximal rigid presentation, with weights $\epsilon_b$, and suppose $\epsilon_a=\epsilon$ is extended-reachable. Its weights form a basis of $\Z^{Q_0}$, by \cite[Section 4.1]{FeiProjection}. Let $c_a$ be the dual row vector, so
\[ \pair{\epsilon_b}{c_a}=\delta_{ba}. \]
Put $\gamma=\gamma(\epsilon)$, and choose a general ordinary kernel $K$ of coweight $\gamma B$, simultaneously realizing the generic Hom dimensions for the finitely many $\epsilon_b$.
For $\sigma=\sum_b n_b\epsilon_b$, with $n_b\in\Z_{\geq0}$, canonical orthogonality and Corollary~\ref{cor:compatible-zero} give $\comp{\epsilon}{\sigma}=-n_a$. The presentation of $\sigma$ is rigid. By \cite[Theorem 3.6]{FeiTrop}, and additivity on its canonical direct sum,
\[ f_K(\sigma)=\hom(\sigma,K) =\hom(\sigma,\gamma B) =\pair{\gamma-c_a}{\sigma}. \]
Homogeneity and continuity extend this equality to the real cone generated by the $\epsilon_b$. Since this cone has nonempty interior and
$f_K(\sigma)=\max_{V\subseteq K}\pair{\dv V}{\sigma}$,
there is a submodule of $K$ with dimension vector $\gamma-c_a$. Indeed, there are only finitely many submodule dimension vectors. Choose an interior point away from the finitely many walls where two of their linear forms agree. A maximizing submodule then remains maximizing in a neighbourhood, and comparison with the linear function $\pair{\gamma-c_a}{\sigma}$ forces its dimension vector to be $\gamma-c_a$. Consequently,
\begin{equation} \label{eq:completion-comparison} \mathbf 0\leq\gamma(\epsilon)-c_a\leq\dv(\gamma(\epsilon)B). \end{equation}

Let $\eta$ be the weight of $M$. For the negative Bongartz completion of a rigid presentation of $\eta$, the distinguished dual vector is $c_a=\dv\rho_+(M)$, by \cite[Theorem 4.11 and Proposition 4.13(2)]{FeiProjection}. The first inequality in \eqref{eq:completion-comparison} proves the positive Schur-reduction bound.

For the negative Schur-reduction bound, first consider $M=I_i=D(Je_i)$. Endomorphisms of $Je_i$ are right multiplications by $e_iJe_i$, and their radical consists of right multiplications by $\mathfrak r_i$. Dualizing and taking the intersection of their kernels gives
\[ \rho_-(I_i)(j)\simeq D\!\left(e_jJe_i/(e_jJe_i)\mathfrak r_i\right),
 \qquad \dim\rho_-(I_i)(j)=\kappa_i(P_j). \]
By Proposition~\ref{prop:projective-denominator} and the local radical bound,
\[ \gamma(I_i)(j)=q_i(P_j)\geq\kappa_i(P_j)=\dim\rho_-(I_i)(j). \]
Thus the injective case does not require applying a completion identity to the negative presentation $\tau_y^{-1}\eta=-\mathbf e_i$.

If $M$ is not injective, put $\epsilon=\tau_y^{-1}\eta$. It has the ordinary rigid representative $\tau^{-1}M$. Take its positive Bongartz completion. The distinguished dual vector is $c_a=-\dv\rho_-(M)$, by \cite[Theorem 4.11 and Proposition 4.13(2)]{FeiProjection}, since $\tau\epsilon=\eta$. The second inequality in \eqref{eq:completion-comparison} gives
\[ \dv\rho_-(M) \leq\dv(\gamma(\epsilon)B)-\gamma(\epsilon) =\gamma(\eta), \]
where the last equality is Theorem~\ref{thm:intrinsic-vector}. The upper bound in \eqref{eq:Schur-comparison} is Corollary~\ref{cor:gamma-positive}. If $\End_J(M)=\kk$, both Schur reductions equal $M$, proving the final assertion.
\end{proof}

The argument uses the dual vectors of Bongartz completions, not the mutation formulas for Schur reductions. The other summands of the completion need not be extended-reachable, because the compatibility theorem permits arbitrary weights.

\begin{proposition}\label{prop:parallel-comparisons}
Let $M$ be an ordinary reachable indecomposable of weight $\eta$, let $z=X(\eta)$, and let $t$ be the current seed. Set
\[ \boldsymbol\kappa(M)=(\kappa_j(M))_{j\in Q_0},\qquad \boldsymbol h(M)=(h_j^J(M))_{j\in Q_0}. \]
Then
\begin{align*}
 \den_t(z)&=\dv M-\boldsymbol h(M),&
 \boldsymbol\kappa(M)&\leq\den_t(z)\leq\dv M,\\
 \gamma(\eta)&=\dv M-\Delta(\eta),&
 \max\{\dv\rho_+(M),\dv\rho_-(M)\}&\leq\gamma(\eta)\leq\dv M.
\end{align*}
All inequalities and the maximum are coordinatewise. The ordered degrees and the corresponding Hessian defects are
\[
\begin{aligned}
 \den_t(z)(j)&=\comp{-\mathbf e_j}{\eta}=d(x_{j;t},z),& h_j(M)&=h(-\mathbf e_j,\eta),\\
 \gamma(\eta)(j)&=\comp{\eta}{-\mathbf e_j}=d(z,x_{j;t}),& \Delta(\eta)(j)&=h(\eta,-\mathbf e_j).
\end{aligned}
\]
\end{proposition}
\begin{proof}
Reachability gives the denominator identity by Theorem~\ref{thm:denominator}; the negative decoration is zero. Theorem~\ref{thm:local-radical-bound} gives
$\dim M(j)-h_j(M)\geq\kappa_j(M)$.
The Hessian identity and its bounds are Theorems~\ref{thm:intrinsic-vector} and~\ref{thm:Schur-comparison}. The ordered-degree assertions follow from Proposition~\ref{prop:compatibility} and Theorem~\ref{thm:ordinary-vector}. For $h(-\mathbf e_j,\eta)$, we choose the identity sequence as normalization, obtaining $h_j(M)$; the other defect identity is part of Theorem~\ref{thm:pair-identity}.
\end{proof}

The same local defect occurs in both formulas. For the denominator vector, it is evaluated on $M$ at the current vertices. For the Hessian vector, choose $w_y\eta=-\mathbf e_i$ and put $\N_j=w(S_j^-)=(N_j,N_j^-)$. Then $\Delta(\eta)(j)=h_i^{J^w}(N_j)$. The comparison can be summarized as follows, under the hypotheses of Proposition~\ref{prop:parallel-comparisons}.
\begin{center}
\renewcommand{\arraystretch}{1.25}
\begin{tabular}{@{}p{3.0cm}p{5.6cm}p{5.6cm}@{}}
\toprule
 & Denominator vector & Hessian vector\\
\midrule
$j$-th degree & $d(x_{j;t},z)$ & $d(z,x_{j;t})$\\
Defect subtracted from $\dv M$ & $h_j(M)$ & $h_i^{J^w}(N_j)$\\
Lower bound & $\boldsymbol\kappa(M)$ & Both $\dv\rho_+(M)$ and $\dv\rho_-(M)$\\
Upper bound & $\dv M$ & $\dv M$\\
\bottomrule
\end{tabular}
\end{center}
No ordering between the denominator and Hessian vectors holds uniformly over all seeds and cluster variables. Such an ordering, applied with the variables interchanged, would force symmetry of the ordered denominator degrees; Example~\ref{ex:equal-denominators} gives an unequal pair. Equivalently,
\[ \den_t(z)-\gamma(\eta)=\Delta(\eta)-\boldsymbol h(M), \]
and the two nonnegative defects need not be ordered. More generally, for distinct extended-reachable indecomposables $\eta,\delta$, subtracting the two instances of \eqref{eq:pair-defect} gives
\[ \comp{\eta}{\delta}-\comp{\delta}{\eta}=h(\delta,\eta)-h(\eta,\delta). \]
The symmetric extension terms cancel, and neither class is a canonical summand of the other. Thus the asymmetry of the compatibility degree is exactly the asymmetry of its defect.

\section{Further directions}\label{sec:outlook}

\subsection{Imaginary presentation classes}
The local defect and the local radical bound apply to all principal components. The construction of $\gamma$, however, requires a negative normalization. An imaginary indecomposable presentation admits none, since the self-extension number of an individual presentation is invariant under extended mutation, whereas negative presentations are rigid. 

\begin{problem}\label{prob:imaginary}
Define a Hessian vector for imaginary presentation weights, prove its ordinary and Auslander--Reiten mutation rules, and determine the cluster-theoretic meaning of the resulting compatibility degree.
\end{problem}
The numerical pairing $\mathcal I_B$ already accepts every integral vector. What remains is to assign the vector to an imaginary weight. The present multiplicity theorem allows imaginary summands in its second argument, but it counts only a specified extended-reachable real summand.

\subsection{The complete canonical decomposition}
Corollary~\ref{cor:acyclic-multiplicity} concerns one real Schur root. Derksen--Weyman \cite{DW} give an effective algorithm for the complete canonical decomposition of a dimension vector of an acyclic quiver. For a finite-dimensional algebra, the Derksen--Fei splitting criterion also gives a finite exhaustive search in principle. The problem is to obtain a structural procedure for presentations that avoids such exhaustive enumeration.

\begin{problem}\label{prob:canonical-algorithm}
Find an effective procedure, in a style of \cite{DW}, that computes the complete canonical decomposition of a presentation weight over a finite-dimensional Jacobian algebra without exhaustive search over all decompositions. 
\end{problem}

\subsection{Module realizations}
Section~\ref{sec:projective} gives a quotient of each projective and a submodule of each injective whose dimension is its Hessian vector. The construction does not establish that the isomorphism class of the realizing module is independent of the deformation direction. In general, the construction in Section~\ref{sec:vectors} gives $n$ local rank comparisons over the Jacobian algebra obtained from a chosen normalization. It does not by itself provide a canonical module over the original algebra of dimension $\gamma(\eta)$.

\begin{problem}
Determine when the Hessian vector is realized by a canonical module or presentation over the original Jacobian algebra. In the projective and injective constructions, determine whether the generic Hom dimensions or Newton polytope are independent of the deformation direction.
\end{problem}

\subsection{Skew-symmetrizable setting}
For a skew-symmetrizable cluster pattern, the ordered denominator degree
$d(x,z)$ is well-defined by \cite[Theorem 6.3]{CL}. Hence, for a cluster
variable $z$ and a seed $t$, one may define
$\gamma_t^d(z)=\bigl(d(z,x_{j;t})\bigr)_j$.
In the standard column convention, mutation at $k$ gives
\[\gamma_{t'}^d(z)(k)=
\max\left\{
\sum_j[b_{jk}]_+\gamma_t^d(z)(j),
\sum_j[-b_{jk}]_+\gamma_t^d(z)(j)
\right\}
-\gamma_t^d(z)(k),
\]
with the other coordinates unchanged. In any initial seed $t_0$ containing
$z=x_{i;t_0}$, the initial value is $\gamma_{t_0}^d(z)=-\mathbf e_i$;
this is \cite[(7.6)--(7.7)]{FZ4}. Thus the cluster-combinatorial part of
the Hessian-vector construction extends directly to the
skew-symmetrizable setting. In particular, if $\mathsf D(t,s)$ denotes
the denominator matrix from $t$ to $s$, then the corresponding reversed
matrix satisfies $\mathsf H^d(t,s)=\mathsf D(s,t)^{\tr}$.

What is missing is a representation-theoretic realization of this
vector by Hessian rank defects. The arguments of this paper use
Jacobian algebras of quivers with potentials, general projective
presentations, and their mutation and Auslander--Reiten translation.
An extension to the skew-symmetrizable setting would therefore require
an analogue of this framework, together with compatible notions of local Hessian defect and
general rank.




\begin{thebibliography}{GLFS20}
\enlargethispage{3\baselineskip}
\bibitem[BFZ05]{BFZ}
A.~Berenstein, S.~Fomin and A.~Zelevinsky, \emph{Cluster algebras III: Upper bounds and double Bruhat cells}, Duke Math. J. \textbf{126} (2005), no.~1, 1--52. \arxiv{math/0305434v3}.

\bibitem[BMR09]{BMR09}
A.~B.~Buan, B.~R.~Marsh and I.~Reiten, \emph{Denominators of cluster variables}, J. Lond. Math. Soc. (2) \textbf{79} (2009), no.~3, 589--611. \arxiv{0710.4335}.

\bibitem[CK06]{CK}
P.~Caldero and B.~Keller, \emph{From triangulated categories to cluster algebras II}, Ann. Sci. \'Ecole Norm. Sup. (4) \textbf{39} (2006), no.~6, 983--1009. \arxiv{math/0510251}.

\bibitem[CL18]{CL18}
P.~Cao and F.~Li, \emph{Positivity of denominator vectors of skew-symmetric cluster algebras}, J. Algebra \textbf{515} (2018), 448--455. \arxiv{1712.06975}.

\bibitem[CL20]{CL}
P.~Cao and F.~Li, \emph{The enough $g$-pairs property and denominator vectors of cluster algebras}, Math. Ann. \textbf{377} (2020), nos.~3--4, 1547--1572. \arxiv{1803.05281}.

\bibitem[CP15]{CP15}
C.~Ceballos and V.~Pilaud, \emph{Denominator vectors and compatibility degrees in cluster algebras of finite type}, Trans. Amer. Math. Soc. \textbf{367} (2015), no.~2, 1421--1439. \arxiv{1302.1052}.

\bibitem[DF15]{DF}
H.~Derksen and J.~Fei, \emph{General presentations of algebras}, Adv. Math. \textbf{278} (2015), 210--237. \arxiv{0911.4913v3}.

\bibitem[DW02]{DW}
H.~Derksen and J.~Weyman, \emph{On the canonical decomposition of quiver representations}, Compos. Math. \textbf{133} (2002), no.~3, 245--265.

\bibitem[DWZ08]{DWZ1}
H.~Derksen, J.~Weyman and A.~Zelevinsky, \emph{Quivers with potentials and their representations I: Mutations}, Selecta Math. (N.S.) \textbf{14} (2008), no.~1, 59--119. \arxiv{0704.0649v4}.

\bibitem[DWZ10]{DWZ2}
H.~Derksen, J.~Weyman and A.~Zelevinsky, \emph{Quivers with potentials and their representations II: Applications to cluster algebras}, J. Amer. Math. Soc. \textbf{23} (2010), no.~3, 749--790. \arxiv{0904.0676v2}.

\bibitem[Fei23a]{FeiComb}
J.~Fei, \emph{Combinatorics of $F$-polynomials}, Int. Math. Res. Not. IMRN (2023), no.~9, 7578--7615. \arxiv{1909.10151v3}.

\bibitem[Fei23b]{FeiTrop}
J.~Fei, \emph{Tropical $F$-polynomials and general presentations}, J. Lond. Math. Soc. (2) \textbf{107} (2023), no.~6, 2079--2120. \arxiv{1911.10513v6}.

\bibitem[Fei25]{FeiRank}
J.~Fei, \emph{On the general ranks of QP representations}, Algebr. Represent. Theory \textbf{28} (2025), no.~1, 47--79. \arxiv{2303.10591v5}.

\bibitem[Fei26]{FeiSyzygy}
J.~Fei, \emph{Syzygy transformations of cluster characters and geometric twists}, \arxiv{2609.10090v1} (2026).

\bibitem[FeiOP]{FeiProjection}
J.~Fei, \emph{On the orthogonal projections}, preprint, 2025; revised 2026. \arxiv{2510.01615v3}.

\bibitem[FG24]{FuGyoda}
C.~Fu and Y.~Gyoda, \emph{Compatibility degree of cluster complexes}, Ann. Inst. Fourier (Grenoble) \textbf{74} (2024), no.~2, 663--718. \arxiv{1911.07193v3}.

\bibitem[FK10]{FK10}
C.~Fu and B.~Keller, \emph{On cluster algebras with coefficients and $2$-Calabi--Yau categories}, Trans. Amer. Math. Soc. \textbf{362} (2010), no.~2, 859--895. \arxiv{0710.3152}.

\bibitem[FZ02]{FZ1}
S.~Fomin and A.~Zelevinsky, \emph{Cluster algebras I: Foundations}, J. Amer. Math. Soc. \textbf{15} (2002), no.~2, 497--529. \arxiv{math/0104151}.

\bibitem[FZ03a]{FZ031}
S.~Fomin and A.~Zelevinsky, \emph{Cluster algebras: Notes for the CDM-03 conference}, Current Developments in Mathematics, 2003, Int. Press, Somerville, MA, 1--34. \arxiv{math/0311493}.

\bibitem[FZ03b]{FZ032}
S.~Fomin and A.~Zelevinsky, \emph{$Y$-systems and generalized associahedra}, Ann. of Math. (2) \textbf{158} (2003), no.~3, 977--1018. \arxiv{hep-th/0111053}.

\bibitem[FZ07]{FZ4}
S.~Fomin and A.~Zelevinsky, \emph{Cluster algebras IV: Coefficients}, Compos. Math. \textbf{143} (2007), no.~1, 112--164. \arxiv{math/0602259}.

\bibitem[GLFS20]{GLFS}
C.~Gei\ss, D.~Labardini-Fragoso and J.~Schr\"oer, \emph{Generic Caldero--Chapoton functions with coefficients and applications to surface cluster algebras}, preprint, 2020. \arxiv{2007.05483v1}.

\bibitem[Kac82]{Kac}
V.~G.~Kac, \emph{Infinite root systems, representations of graphs and invariant theory, II}, J. Algebra \textbf{78} (1982), no.~1, 141--162.

\bibitem[LS15]{LS}
K.~Lee and R.~Schiffler, \emph{Positivity for cluster algebras}, Ann. of Math. (2) \textbf{182} (2015), no.~1, 73--125. \arxiv{1306.2415}.

\bibitem[Rea07]{Rea07}
N.~Reading, \emph{Clusters, Coxeter-sortable elements and noncrossing partitions}, Trans. Amer. Math. Soc. \textbf{359} (2007), no.~12, 5931--5958. \arxiv{math/0507186}.

\bibitem[RS16]{RSpeyer}
N.~Reading and D.~E.~Speyer, \emph{Combinatorial frameworks for cluster algebras}, Int. Math. Res. Not. IMRN (2016), no.~1, 109--173. \arxiv{1111.2652}.

\bibitem[RS18]{RS}
N.~Reading and S.~Stella, \emph{Initial-seed recursions and dualities for $d$-vectors}, Pacific J. Math. \textbf{293} (2018), no.~1, 179--206. \arxiv{1409.4723}.

\bibitem[Sch92]{Sch92}
A.~Schofield, \emph{General representations of quivers}, Proc. London Math. Soc. (3) \textbf{65} (1992), no.~1, 46--64.

\bibitem[Sta]{Stacks}
The Stacks Project Authors, \emph{The Stacks Project}, Lemma 37.24.3, \href{https://stacks.math.columbia.edu/tag/054X}{Tag 054X}.

\bibitem[Yur24]{Yur24}
T.~Yurikusa, \emph{Denominator vectors and dimension vectors from triangulated surfaces}, J. Algebra \textbf{641} (2024), 620--647. \arxiv{2305.08097}.

\end{thebibliography}
\end{document}